\documentclass[aop,MSNbibl,dvips]{arximspdf}
\usepackage{mathbh}
\usepackage{graphicx}

\doi{10.1214/13-AOP890} 
\volume{42}
\issue{5}
\pubyear{2014}
\firstpage{1769}
\lastpage{1808}

\makeatletter
\newtheorem{theorem}{Theorem}
\newproclaim{remark}[theorem]{Remark}
\newtheorem{lemma}[theorem]{Lemma}
\newtheorem{proposition}[theorem]{Proposition}
\newproclaim{definition}[theorem]{Definition}
\newtheorem{conjecture}[theorem]{Conjecture}
\def\one{\mathbh{1}}
\def\widebar{\overline}
\makeatother

\begin{document}
\begin{frontmatter}

\title{Critical Gaussian multiplicative chaos: Convergence~of the
derivative martingale}
\runtitle{Critical Gaussian multiplicative chaos}

\begin{aug}
\author[A]{\fnms{Bertrand} \snm{Duplantier}\ead[label=e1]{bertrand.duplantier@cea.fr}\thanksref{t1}},
\author[B]{\fnms{R\'{e}mi} \snm{Rhodes}\ead[label=e2]{rhodes@ceremade.dauphine.fr}\thanksref{t2}},\\
\author[C]{\fnms{Scott} \snm{Sheffield}\ead[label=e3]{sheffield@math.mit.edu}\thanksref{t3}}
\and
\author[B]{\fnms{Vincent} \snm{Vargas}\corref{}\ead[label=e4]{vargas@ceremade.dauphine.fr}\thanksref{t4}}
\runauthor{Duplantier, Rhodes, Sheffield and Vargas}
\affiliation{Institut de Physique Th\'{e}orique, Universit\'e
Paris-Dauphine,\\ Massachusetts Institute of Technology and Universit\'e Paris-Dauphine}
\address[A]{B. Duplantier\\
Institut de Physique Th\'{e}orique\\
CEA/Saclay\\
F-91191 Gif-sur-Yvette Cedex\\
France\\
\printead{e1}} 
\address[B]{R. Rhodes\\
V. Vargas\\
Universit\'e Paris-Dauphine\\
Ceremade, UMR 7564\\
Place de Mar\'{e}chal de Lattre de Tassigny\\
75775 Paris Cedex 16\\
France\\
\printead{e2}\\
\phantom{E-mail:\ }\printead*{e4}}
\address[C]{S. Sheffield\\
Department of Mathematics\\
Massachusetts Institute of Technology\\
Cambridge, Massachusetts 02139\\
USA\\
\printead{e3}}
\end{aug}
\thankstext{t1}{Supported in part by Grant ANR-08-BLAN-0311-CSD5 and by
the MISTI MIT-France Seed Fund.}
\thankstext{t2}{Supported in part by Grant ANR-11-JCJC CHAMU.}
\thankstext{t3}{Supported in part by NSF Grants DMS-06-4558 and OISE
0730136 and by the MISTI MIT-France Seed Fund.}
\thankstext{t4}{Supported in part by Grant ANR-11-JCJC CHAMU.}

\received{\smonth{7} \syear{2012}}
\revised{\smonth{9} \syear{2013}}

%
\begin{abstract}
In this paper, we study Gaussian multiplicative chaos in the critical
case. We show that the so-called derivative martingale, introduced in
the context of branching Brownian motions and branching random walks,
converges almost surely (in all dimensions) to a random measure with
full support. We also show that the limiting measure has no atom. In
connection with the derivative martingale, we write explicit
conjectures about the glassy phase of log-correlated Gaussian
potentials and the relation with the asymptotic expansion of the
maximum of log-correlated Gaussian random variables.
\end{abstract}

%
\pdfkeywords{60G57, 60G15, 60D05, Gaussian multiplicative chaos, Liouville quantum gravity, maximum of log-correlated fields}

\end{frontmatter}

\section{Introduction}\label{sec1}
\subsection{Overview}\label{sec1.1}
In the 1980s, Kahane \cite{cfKah} developed a continuous parameter
theory of multifractal random measures, called Gaussian multiplicative
chaos; this theory emerged from the need to define rigorously the limit
lognormal model introduced by Mandelbrot \cite{cfMan} in the context
of turbulence. His efforts were followed by several authors \cite
{allez,Bar,bacry,Fan,sohier,rhovar,cfRhoVar} coming up with various
\mbox{generalizations} at different scales. This family of random fields has
found many applications in various fields of science, especially in
turbulence and in mathematical finance. Recently, the authors in \cite
{cfDuSh} constructed a probabilistic and geometrical framework for
Liouville quantum gravity and the so-called
Knizhnik--Polyakov--Zamolodchikov (KPZ) equation \cite{cfKPZ}, based
on the two-dimensional Gaussian free field (GFF); see \cite
{cfDa,DFGZ,DistKa,cfDuSh,GM,cfKPZ,Nak} and references therein. In
this context, the KPZ formula has been proved rigorously \cite
{cfDuSh}, as well as in the general context of Gaussian multiplicative
chaos \cite{cfRhoVar}; see also \cite{Benj} in the context of
Mandelbrot's multiplicative cascades. This was done in the standard
case of Liouville quantum gravity, namely strictly below the critical
value of the GFF coupling constant $\gamma$ in the Liouville conformal
factor, that is, for $\gamma<2$ (in a chosen normalization). Beyond
this threshold, the standard construction yields vanishing random
measures \cite{PRL,cfKah}. The issue of mathematically constructing
singular Liouville measures beyond the phase transition (i.e., for
$\gamma>2$) and deriving the corresponding (nonstandard dual) KPZ
formula has been investigated in \cite{BJRV,Duphouches,PRL}, giving
the first mathematical understanding of the so-called \textit{duality} in
Liouville quantum gravity; see \cite
{Al,Amb,Das,BDMan,Dur,Jain,Kleb1,Kleb2,Kleb3,Kostovhouches} for an
account of physical motivations. However, the rigorous construction of
random measures \textit{at criticality}, that is, for $\gamma=2$, does
not seem to ever have been carried out.

As stated above, once the Gaussian randomness is fixed, the standard
Gaussian multiplicative chaos describes a random positive measure for
each $\gamma< 2$ but yields~$0$ when $\gamma=2$. Naively, one might
therefore guess that $-1$ times the \textit{derivative} at $\gamma= 2$
would be a random positive measure. This intuition leads one to
consider the so-called \textit{derivative martingale}, formally
obtained by differentiating the standard measure w.r.t. $\gamma$ at
$\gamma=2$, as explained below. In the case of branching Brownian
motions \cite{neveu}, or of branching random walks \cite{BiKi,Kyp} (see
also \cite{AidShi} for a recent different but equivalent construction),
the construction of such an object has already been carried out
mathematically. In the context of branching random walks, the
derivative martingale was introduced in the study of the fixed points
of the smoothing transform at criticality (the smoothing transform is a
generalization of Mandelbrot's $\star$-equation for discrete
multiplicative cascades; see also \cite{Biggf}). Our construction will
therefore appear as a continuous analogue of those works in the context
of Gaussian multiplicative chaos.

Besides the 2D-Liouville Quantum Gravity framework (and the KPZ
formula), many other important models or questions involve Gaussian
multiplicative chaos of log-correlated Gaussian fields in all
dimensions. Let us mention the glassy phase of log-correlated random
potentials (see \cite{arguin,CarDou,Fyo,rosso}) or the asymptotic
expansion of the maximum of log-correlated random variables; see \cite
{bramson,DingZei}. In all these problems, one of the key tools is the
derivative martingale at the critical point $\gamma^2=2d$ (where $d$ is
the dimension), whose construction is precisely the purpose of this paper.



In dimension $d$, a standard Gaussian multiplicative chaos is a random
measure that can be written formally, for any Borelian set $A\subset
\mathbb R^d$, as
%
\begin{equation}
\label{measintro} M^\gamma(A)=\int_Ae^{\gamma X(x)-(\gamma^2/2) \mathbb E [X^2(x)]}\,dx,
\end{equation}
where $X$ is a centered log-correlated Gaussian field
\[
\mathbb E \bigl[X(x)X(y)\bigr]=\ln_+\frac{1}{|x-y|}+g(x,y)
\]
with $\ln_+(x)=\max(\ln x, 0)$ and $g$ a continuous bounded function
over $\mathbb{R}^d\times\mathbb{R}^d$. Although such an $X$ cannot
be defined as a
random function (and may be a random \textit{distribution}, like the GFF),
the measures can be rigorously defined all for \mbox{$\gamma^2<2d$} using a
straightforward limiting procedure involving a time-indexed family of
improving approximations to $X$ \cite{cfKah}, as we will review in
Section~\ref{setup}. By contrast, it is well known that for $\gamma
^2\geq2d$ the measures constructed by this procedure are identically
zero \cite{cfKah}. Other techniques are thus required to create
similar measures beyond the critical value $\gamma^2=2d$ \cite
{BJRV,Duphouches,PRL}.

Roughly speaking, the derivative martingale 
is defined as (recall that $\gamma=\sqrt{2d}$ is the critical value)
%
\begin{eqnarray}\label{derivintro}
M'(A) &:=& -\frac{\partial}{\partial\gamma} \bigl[
M^\gamma(A) \bigr]_{\gamma= \sqrt{2d}}
\nonumber\\[-8pt]\\[-8pt]
& =& \biggl[ \int_A
\bigl(\gamma \mathbb E \bigl[X^2(x)\bigr]-X(x) \bigr)e^{\gamma X(x)-(\gamma^2/2) \mathbb E
[X^2(x)]}
\,dx \biggr]_{\gamma= \sqrt{2d}}.\nonumber
\end{eqnarray}
Here we have differentiated the measure $M^\gamma$ in~(\ref{measintro})
with respect to the parameter~$\gamma$ to obtain the above expression
(\ref{derivintro}). Note that this is the same as~(\ref{measintro})
except for the factor $ (\gamma \mathbb E [X^2(x)]-X(x) )$.
To give the
reader some intuition, we remark that we will ultimately see that the
main contributions to $M'(A)$ come from locations $x$ where this factor
is positive but relatively close to zero (on the order of $\sqrt {\mathbb E[X^2(x)]}$) which correspond to locations $x$ where $X(x)$ is
nearly maximal. Indeed, in what follows, the reader may occasionally
wish to forget the derivative interpretation of~(\ref{derivintro}) and
simply view $ (\gamma \mathbb E [X^2(x)]-X(x) )$ as the
factor by which
one rescales~(\ref{measintro}) in order to ensure that one obtains a
nontrivial measure (instead of zero) when using the standard limiting
procedure.

\begin{figure}

\includegraphics{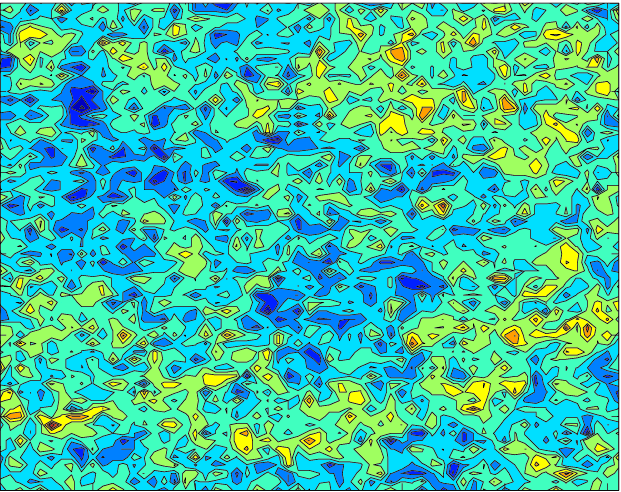}

\caption{Height landscape of the derivative martingale measure plotted
with a logarithmic scale color-bar, showing that the measure is very
``peaked'' (for $t=12$, a multiplicative factor of about $10^{8}$
stands between extreme values, i.e., between warm and cold colors).}
\label{fig1}
\end{figure}

In a sense, the measures $M^\gamma$ in~(\ref{measintro}) become more
concentrated as $\gamma^2$ approaches~$2d$. (They assign full measure
to a set of Hausdorff dimension $d-\gamma^2/2$, which tends to zero as
$\gamma^2 \to2d$.) It is therefore natural to wonder how concentrated
the $\gamma^2 = 2d$ measure will be (see Figure~\ref{fig1} for a simulation of the landscape). In particular, it is natural to
wonder whether it possesses \textit{atoms} (in which case it could in
principle assign full measure to a countable set). In our context, we
will answer \textit{in the negative}. At the time we posted the first
version of this manuscript online, this question was open in the
context of discrete models as well as continuous models. However, a
proof of the nonatomicity of the discrete cascade measures was posted
very shortly afterward in \cite{BKNSW}, which uses a
method independent of our proof. Since our proof is based on a
continuous version of the spine decomposition, as developed in the
context of branching random walks, we expect that it can be adapted to
these other models as well.

Roughly speaking, the reason that establishing nonatomicity in
critical models is nontrivial is that proofs of nonatomicity for
(noncritical) multiplicative chaos usually rely on the existence of
moments higher than $1$ (see \cite{daley}) and the scaling relations of
multifractal random measures; see, for example, \cite{allez}. At
criticality, the random measures involved (cascades, branching random
walks, or Gaussian multiplicative chaos) no longer possess finite
moments of order $1$, and the scaling relations become useless.

To explain this issue in more detail, we recall that it is proved in
\cite{daley} that a stationary random measure $M$ over $\mathbb{R}^d$
is almost
surely nonatomic if ($C$ stands here for the unit cube of $\mathbb{R}^d$)
%
\begin{equation}
\label{critdaley} \forall\delta>0\qquad n^d \mathbb P \bigl(M
\bigl({n}^{-1}C\bigr)>\delta \bigr)\to0\qquad \mbox{as }n\to\infty.
\end{equation}
When $M=M_\gamma$ for $0<\gamma^2<2d$, a computable property of
$M_\gamma$ is its power-law spectrum $\xi$ characterized by
%
\begin{equation}
\label{critsp} \mathbb E \bigl[ \bigl(M^\gamma\bigl({n}^{-1}C
\bigr) \bigr)^q \bigr]\simeq D_q {n^{-\xi
(q)}}\qquad\mbox{as }n\to\infty
\end{equation}
for all those $q$ making the above expectation finite, that is, $q\in
[0,\frac{2d}{\gamma^2}[$. It matches
%
\begin{equation}
\label{critpls} \xi(q)= \biggl(d+\frac{\gamma^2}{2} \biggr)q-\frac{\gamma^2}{2}q^2.
\end{equation}
Using the Markov inequality in~(\ref{critdaley}),~(\ref{critsp})
obviously yields for $q\in[0,\frac{2d}{\gamma^2}[$
\[
n^d \mathbb P \bigl(M^\gamma\bigl({n}^{-1}C
\bigr)>\delta \bigr)\lesssim\frac
{D_q}{\delta
^q}n^{d-\xi(q)}.
\]
Therefore, the nonatomicity of the measure boils down to finding a $q$
such that the power-law spectrum is strictly larger than $d$:
\begin{itemize}
\item In the subcritical situation $\gamma^2< 2d$, the function $\xi$
increases on $ [0,1]$ from $0$ to $d$. Such a $q$ is necessarily larger
than $1$, and a straightforward computation shows that any $q\in\,]1,
\frac{2d}{\gamma^2}[$ suffices.
\item For $\gamma^2= 2d$, relations~(\ref{critsp}) and (\ref
{critpls}) should remain valid only for $q<1$. Therefore, the
subcritical strategy fails because the power-law spectrum achieves its
maximum $d$ at $q=1$. It is tempting to try to replace the gauge
function $x\mapsto x^q$ by something that could be more appropriate at
criticality like $x\mapsto x\ln(1+x)^q$, etc. However, the fact that
the measure does not possess a moment of order $1$ (see Proposition
\ref
{moments} below) shows that there is no way of changing the gauge so as
to make $\xi$ go beyond $d$.
\end{itemize}
More sophisticated machinery is thus necessary to investigate
nonatomicity at criticality. Indeed, we expect the derivative
martingale to assign full measure to a (random) Hausdorff set of
dimension $0$, indicating that the measure is in some sense just
``barely'' nonatomic.


Let us finally mention the interesting work of \cite{tecu} where the
author constructs on the unit circle ($d=1$) a classical Gaussian
multiplicative Chaos given by the exponential of a field $X$ such that
for each $\varepsilon$ the covariance of $X$ at points $x$ and $y$ lies
strictly between $(2-\varepsilon) \ln_+\frac{1}{|x-y|} $ and $2\ln
_+\frac{1}{|x-y|}$ when $|x-y|$ is sufficiently small. In some sense,
his construction is a \textit{near} critical construction, different from
the measures constructed here. This is illustrated by the fact that the
measures in~\cite{tecu} possess moments of order $1$ (and even belong
to $L \log L$), which is atypical for the critical multiplicative chaos
associated to log-correlated random variables.

In this paper, we tackle the problem of constructing random measures at
criticality for a large class of log-correlated Gaussian fields in any
dimension, the covariance kernels of which are called $\star$-scale
invariant kernels.
This approach allows us to link the measures under consideration to a
functional equation, the $\star$-equation, giving rise to several
conjectures about the glassy phase of log-correlated Gaussian
potentials and about the three-terms expansion of the maximum of
log-correlated Gaussian variables.

Another important family of random measures is the class defined by
taking~$X$ to be the Gaussian Free Field (GFF) with free or Dirichlet
boundary conditions on a planar domain, as in \cite{cfDuSh}; see also
\cite{She07} for an introduction to the GFF. The measures defined in
this way are also known as the (critical) Liouville quantum gravity
measures, and are closely related to conformal field theory, as well as
various \mbox{2-}dimensional discrete random surface models and their
scaling limits. Although the Gaussian free field is in some sense a
log-correlated random field, it does not fall exactly into the
framework of this paper, which deals with translation invariant random
measures (defined on all of $\mathbb R^2$ or $\mathbb R^d$) that can be
approximated in a particular way (via the $\star$-equation). Although
some of the arguments of this paper can be easily extended to settings
where the strict translation invariance requirement for $X$ is relaxed
(e.g., $X$ is the Gaussian free field on a disk), we will still need
additional arguments to show that the derivative martingale associates
a unique nonatomic random positive measure to a given instance of the
GFF almost surely, that this measure is independent of the particular
approximation scheme used, and that this measure transforms under
conformal maps in the same way as the $\gamma<2$ measures constructed
in \cite{cfDuSh}. For the sake of pedagogy, this other part of our
work will appear in a companion paper. For the time being, we just
announce that all the results of this paper are valid for the GFF construction.

\subsection{Physics literature: History and motivation}\label{sec1.2}
It is interesting to pause for a moment and consider the physics
literature on Liouville quantum gravity. We first remark that the
noncritical case, with $d=2$ and $\gamma< 2$, was treated in \cite
{cfDuSh}, which contains an extensive overview of the physics
literature and an explanation of the relationships (some proved, some
conjectural) between random measures and discrete and continuum random
surfaces. Roughly speaking, when one takes a random two-dimensional
manifold and conformally maps it to a disk, the image of the area
measure is a random measure on the disk that should correspond to an
exponential of a log-correlated Gaussian random variable (some form of
the GFF). From this point of view, many of the physics results about
discrete and continuum random surfaces can be interpreted as
predictions about the behavior of these random measures, where the
value of $\gamma< 2$ depends on the particular physical model in question.

There is also a physics literature focusing on the critical case
$\gamma
=2$, which we expect to be related to the measure constructed in this
paper. This section contains a brief overview of the results from this
literature, as appearing in, for example, \cite{BKZ,GM,GZ,GrossKleban,GrossM,GubserKleban,KKK,Kleb2,kostov91,kostov92,KS,Parisi,Polch,sugino}.
Most of the results surveyed in this section have not yet been
established or understood in a mathematical sense.

The critical case $\gamma=2$ corresponds to the value $c=1$ of the
so-called \textit{central charge} $c$ of the conformal field theory
coupled to gravity, via the famous KPZ result \cite{cfKPZ},
\[
\gamma=\frac{1}{\sqrt{6}}(\sqrt{25-c}-\sqrt{1-c}).
\]
Discrete critical statistical physical models having $c=1$ then include
one-dimensional matrix models [also called ``matrix quantum mechanics''
(MQM)] \cite
{BKZ,GM,GZ,GrossKleban,GubserKleban,KKK,Kleb2,Parisi,Polch,sugino}, the
so-called $O(n)$ loop model on a random planar lattice for $n=2$ \cite
{kostov91,kostov92,Kostovhouches,KS} and the $Q$-state Potts model on
a random lattice for $Q=4$ \cite{BBM,Daul,Ey}. For an introduction to
the above mentioned 2D statistical models, see, for example, \cite{Nienhuis}.

In the continuum, a natural coupling also exists between Liouville
quantum gravity and the Schramm--Loewner evolution SLE$_\kappa$ for
$\gamma=\sqrt{\kappa}$, rigorously established for $\kappa<4$ \cite
{PRL1,She11}. Thus the critical value $\gamma=2$ corresponds to the
special SLE parameter value $\kappa=4$, above which the SLE$_\kappa$
curve no longer is a simple curve, but develops double points at all scales.

The standard $c=1, \gamma=2$ Liouville field theory \cite
{BKZ,GM,GZ,GrossKleban,KKK,Kleb2,Parisi,Polch} involves violations of
scaling by \textit{logarithmic factors}. For example, the partition
function (number) of genus $0$ random surfaces of area $\mathcal A$
grows as \cite{GrossKleban,KKK}
\[
\mathcal Z\propto\exp(\mu\mathcal A) \mathcal A^{-3} (\log\mathcal
A)^{-2},
\]
where $\mu$ is a nonuniversal growth constant depending on the (planar
lattice) regularization. The area exponent ($-$3) is universal for a
$c=1$ central charge, while the subleading logarithmic factor is
attributed to the unusual dependence on the Liouville field $\varphi$
(equivalent to $X$ here) of the so-called ``tachyon field'' $T(\varphi
)\propto\varphi  e^{2\varphi}$ \cite{KKK,Kleb2,Polch}. Its integral
over a ``background'' Borelian set $A$ generates the quantum area
$\mathcal A =\int_A T(\varphi) \,dx $, that we can recognize as the
formal heuristic expression for the \textit{derivative measure} (\ref
{derivintro}) introduced above.

At $c=1$, a proliferation of large ``bubbles'' (the so-called ``baby
universes'' which are relatively large amounts of area cut off by
relatively small bottlenecks) is generally anticipated in the bulk of
the random surface \cite{GrossKleban,Jain,kostov91}, or at its boundary
in the case of a disk topology \cite{kostov92,KS}. We believe that this
should correspond to the fact that the measure we construct is
concentrated on a set of Hausdorff dimension zero.

However, the introduction of higher trace terms \cite
{GubserKleban,Kleb2,sugino} in the action of the $c=1$ matrix model of
two-dimensional quantum gravity is known to generate a ``nonstandard''
random surface model with an even stronger concentration of
bottlenecks. (See also the related detailed study of a MQM model for a
$c=1$ string theory with vortices in \cite{KKK}.) As we shall see
shortly, these nonstandard constructions do not seem to correspond to
our model, at least not so directly.
In these constructions, one encounters a new critical behavior of the
random surface, with a~\textit{critical} proliferation of spherical
bubbles connected one to another by microscopic ``wormholes.'' This is
reminiscent of the construction for $c<1, \gamma<2$ of the \textit
{dual} phase of Liouville quantum gravity \cite
{Al,Amb,Das,Dur,Kleb1,Kleb2,Kleb3}, where the associated random measure
develops \textit{atoms} \cite{BJRV,Duphouches,PRL}.

The partition function of the nonstandard $c=1$ (genus zero) random
surface then scales as a function of the area $\mathcal A$ as \cite
{GubserKleban,KKK,Kleb2,sugino}
\[
\mathcal Z\propto\exp\bigl(\mu'\mathcal A\bigr) \mathcal
A^{-3}
\]
with an apparent suppression of logarithmic terms. This has been
attributed to the appearance for $c=1$ of a tachyon field of the
\textit
{atypical} form $T(\varphi)\propto e^{2\varphi}$ \mbox{\cite{GubserKleban,KKK,Kleb3}}. Heuristically, this would seem to correspond
to a measure of type~(\ref{measintro}), but we know that the latter
vanishes for $\gamma=2$. (See Proposition~\ref{propstand} below.) The
literature about the analogous problem of branching random walks \cite
{AidShi,HuShi} also suggests for $\gamma=2$ a \textit{logarithmically
renormalized} measure obtained by multiplying by $\sqrt{\log(1/\varepsilon)}=\sqrt{t}$ the object [see~(\ref{defchaos}) below]
whose limit is taken in~(\ref{measintro}), but we expect this to
converge (up to constant factor) to the same measure as the derivative
martingale~(\ref{derivintro}). In order to model the nonstandard
theory, it might be necessary to modify the measures introduced here by
explicitly introducing ``atoms'' on top of them, using the procedure
described in \cite{BJRV,Duphouches,PRL} for adding atoms to $\gamma<2$
random measures. In the approach of \cite{BJRV,Duphouches,PRL}, the
``dual Liouville measure'' corresponding to $\gamma< 2$ involves
choosing a Poisson point process from $\eta^{-\alpha-1}\,d\eta M_\gamma
(dx)$, where $\alpha= \gamma^2/4 \in(0,1)$, and letting each point
$(\eta,x)$ in this process indicate an atom of size $\eta$ at location
$x$. When $\gamma=2$ and $\alpha= 1$, we can replace $M_\gamma$ with
the $M'$ of~(\ref{derivintro}) and use the same construction; in this
case (since $\alpha= 1$) the measure a.s. assigns infinite mass to
each positive-Lebesgue-measure $A \in\mathcal B(\mathbb R^d)$.
However, one may use standard L\'evy compensation to produce a random
distribution, assigning a finite value a.s. to each fixed $A \in
\mathcal B(\mathbb R^d)$ with a positive atom of size $\eta$ at
location $x$ corresponding to each $(\eta,x)$ in the Poisson point
process. We suspect that that this construction is somehow equivalent
to the continuum random measure associated with the nonstandard
$c=1,\gamma=2$ Liouville random surface with enhanced bottlenecks, as
described in \cite{GubserKleban,KKK,sugino}.

Finally, we note that the \textit{boundary} critical Liouville quantum
gravity poses similar challenges. A subtle difference in logarithmic
boundary behavior is predicted between the so-called \textit{dilute}
and \textit{dense} phases of the $O(2)$ model on a random disk \cite
{kostov92,KS}, which thus may differ in their boundary bubble
structure. It also remains an open question whether the results about
the conformal welding of two boundary arcs of random surfaces to
produce SLE, as described in \cite{She11}, can be extended to the case
$\gamma=2$.


\section{Setup}\label{setup}
\subsection{Notation}\label{sec2.1}
For a Borelian set $A\subset\mathbb{R}^d$, $\mathcal{B}(A)$ stands
for the
Borelian sigma-algebra on $A$. All the considered fields are
constructed on the same probability space $(\Omega,\mathcal{F},
\mathbb P )$. We
denote by $ \mathbb E $ the corresponding expectation.

\subsection{\texorpdfstring{$\star$}{star}-scale invariant kernels}\label{sec2.2}
Here we introduce the Gaussian fields that we will use throughout the
papers. We consider a family of centered stationary Gaussian processes
$((X_{t}(x))_{x \in\mathbb{R}^d})_{t\geq0}$ where, for each $t\geq
0$, the
process $(X_{t}(x))_{x \in\mathbb{R}^d}$ has covariance given by
%
\begin{equation}
\label{corrX} K_t(x)= \mathbb E \bigl[X_{t}(0)X_{t}(x)
\bigr]= \int_1^{e^t}\frac{k(ux)}{u} \,du
\end{equation}
for some covariance kernel $k$ satisfying $k(0)=1$, of class $C^1$ and
vanishing outside a compact set (actually this latter condition is not
necessary but it simplifies the presentation). The $C^1$ condition is
technical and ensures that for $x \neq y$ we have a nice description
of the joint law of the couple $(X_t(x),X_t(y))_{ t \geq0}$; see Lemma
\ref{decomp} below (this condition could also be relaxed to some
extent). We also assume that the process $(X_t(x)-X_s(x))_{x\in\mathbb{R}^d}$
is independent of the processes $ ((X_u(x))_{x\in\mathbb
{R}^d} )_{u\leq
s}$ for all $s<t$. Put in other words, the mapping $t\mapsto X_t(\cdot
)$ has independent increments. Such a construction of Gaussian
processes is carried out in \cite{allez}.
For $\gamma\geq0$, we consider the approximate Gaussian multiplicative
chaos $M^\gamma_{t}(dx)$ on $\mathbb{R}^d$,
%
\begin{equation}
\label{defchaos} M^\gamma_{t}(dx)=e^{\gamma X_{t}(x)-(\gamma^2/2)E[X_{t}(x)^2]}\,dx.
\end{equation}

It is well known \cite{allez,cfKah} that, almost surely, the family
of random measures $(M^\gamma_t)_{t>0}$ weakly converges as $t\to
\infty
$ toward a random measure $M^\gamma$, which is nontrivial if and only
if $\gamma^2<2d$. The purpose of this paper is to investigate the phase
transition, that is, $\gamma^2=2d$. Recall that we have:
%
\begin{proposition}
For $\gamma^2=2d$, the standard construction~(\ref{defchaos}) yields a
vanishing limiting measure
%
\begin{equation}
\label{eqngamma2limitzero} \lim_{t \to\infty} M^{\sqrt{2d}}
_t(dx)=0\qquad\mbox{almost surely}.
\end{equation}
\end{proposition}

Let us also mention that the authors in \cite{allez} have proved that,
for $\gamma^2<2d$, the measure $M^\gamma$ satisfies the following scale
invariance relation, called $\star$-equation:
%
\begin{definition}[(Log-normal $\star$-scale invariance)]\label{def1}
The random Radon measure $M^\gamma$ is lognormal $\star$-scale
invariant: for all $0<\varepsilon\leq1$, $M^\gamma$ obeys the
cascading rule
%
\begin{eqnarray}\label{star}
&& \bigl(M^\gamma(A) \bigr)_{A\in\mathcal{B}(\mathbb{R}^d)}
\nonumber\\[-4pt]\\[-12pt]
&&\qquad \stackrel{\mathrm{law}}
{=} \biggl(\int_Ae^{\gamma X_{\ln(1/\varepsilon)}(r)-(\gamma^2/2) \mathbb E
[X_{\ln(1/\varepsilon)}(r)^2]}\varepsilon^{d}
M^{\gamma,\varepsilon}(dr) \biggr)_{A\in\mathcal{B}(\mathbb{R}^d)},\nonumber
\end{eqnarray}
where $X_{\ln(1/\varepsilon)}$ is the Gaussian process
introduced in~(\ref{corrX}), and $M^{\gamma,\varepsilon}$ is a random
measure independent from $X_{\ln(1/\varepsilon)}$ satisfying the
scaling relation
%
\begin{equation}
\label{star1} \bigl(M^{\gamma,\varepsilon}(A) \bigr)_{A\in\mathcal{B}(\mathbb
{R}^d)}\stackrel {\mathrm{law}}
{=} \biggl(M^\gamma\biggl(\frac{A}{\varepsilon}\biggr) \biggr)_{A\in\mathcal
{B}(\mathbb{R}^d)}.
\end{equation}
\end{definition}

Intuitively, this relation means that when zooming in the measure $M$,
one should observe the same behavior up to an independent Gaussian
factor. It has some canonical meaning since it is the exact continuous
analog of the smoothing transformation intensively studied in the
context of Mandelbrot's multiplicative \mbox{cascades} \cite{durrett} or
branching random walks \cite{Biggf,Liu}.

Observe that this equation perfectly makes sense for the value $\gamma
^2=2d$. Therefore, to define a natural Gaussian multiplicative chaos at
the value $\gamma^2=2d$, one has to look for a solution to this
equation when $\gamma^2=2d$ and conversely, each random measure
candidate for being a Gaussian multiplicative chaos at the value
$\gamma
^2=2d$ must satisfy this equation.

\begin{remark}
The main motivation for considering $\star$-scale invariant kernels is
the connection between the associated random measures and the $\star
$-equation. Nevertheless, we stress that our proofs can be easily
adapted to other Gaussian\vadjust{\goodbreak} multiplicative chaos associated to
log-correlated Gaussian fields ``\`a la Kahane'' \cite{cfKah}: in
particular, we can construct the derivative martingale associated to
exact scale invariant kernels \cite{bacry,rhovar} or the Gaussian Free
Field in a bounded domain.
\end{remark}
%

\section{Derivative martingale}\label{secder}

One way to construct a solution to the $\star$-equation at the critical
value $\gamma^2=2d$ is to introduce the derivative martingale
$M_{t}'(dx)$ defined by
\begin{eqnarray*}
M_{t}'(dx)&:=& \bigl(\sqrt{2d} t-X_{t}(x)
\bigr)e^{\sqrt{2d}X_{t}(x)-d \mathbb E
[X_{t}(x)^2]}\,dx.
\end{eqnarray*}
It is plain to see that, for each open bounded set $A\subset\mathbb
{R}^d$, the
family $(M_{t}'(A))_t$ is a martingale. Nevertheless, it is not
nonnegative. It is therefore not obvious that such a family converges
toward a (nontrivial) positive limiting random variable. The following
theorem is the main result of this section:

\begin{theorem}\label{mainderiv}
For each bounded open set $A\subset\mathbb{R}^d$, the martingale\break
$(M'_t(A))_{t\geq0}$\vadjust{\goodbreak} converges almost surely toward a positive random
variable denoted by $M'(A)$, such that $M'(A)>0$ almost surely.
Consequently, almost surely, the (locally signed) random measures
$(M'_t(dx))_{t\geq0}$ converge weakly as $t\to\infty$ toward
a~positive random measure $M'(dx)$. This limiting measure has full
support and is atomless. Furthermore, the measure $M'$ is a solution to
the $\star$-equation~(\ref{star}) with $\gamma=\sqrt{2d}$.
\end{theorem}

Since $M_{t}'(dx)$ is not uniformly nonnegative when $t < \infty$,
there are several complications involved in establishing its
convergence to a nonnegative limit (let alone the nontriviality of
the limit). We have to introduce some further tools to study its
convergence. These tools have already been introduced in the context of
discrete multiplicative cascade models in order to study the
corresponding derivative martingale; see \cite{BiKi}.

We denote by $\mathcal{F}_t$ the sigma algebra generated by $\{
X_s(x);s\leq t,x\in\mathbb{R}^d\}$. Given a Borelian set $A\subset
\mathbb{R}^d$ and
parameters $t,\beta>0$, we introduce the random variables
\begin{eqnarray*}
Z^{\beta}_t(A)&=&\int_A\bigl(
\sqrt{2d}t-X_{t}(x)+\beta\bigr)\one_{\{\tau
^\beta>t\}
} e^{\sqrt{2d}X_{t}(x)-d \mathbb E [X_{t}(x)^2]} \,dx,
\\
\widetilde{Z}^{\beta}_t(A)&=&\int_A
\bigl(\sqrt{2d}t-X_{t}(x)\bigr)\one_{\{
\tau
^\beta>t\}} e^{\sqrt{2d}X_{t}(x)-d \mathbb E [X_{t}(x)^2]}
\,dx,
\end{eqnarray*}
where, for each $x\in A$, $\tau^\beta(x)$ is the $(\mathcal
{F}_t)_t$-stopping time defined by
\[
\tau^\beta(x)=\inf\bigl\{u>0,X_u(x)-\sqrt{2d}u>\beta\bigr
\}.
\]
In the sequel, when the context is clear, we will drop the $x$
dependence in $\tau^\beta(x)$. What is the relation between $Z^{\beta
}_t(A)$ and $M_{t}'(A)$? Roughly speaking, we will show that the
convergence of $M'_t(A)$ as $t\to\infty$ toward a nontrivial object
boils down to proving the convergence of $Z^{\beta}_t(A)$ toward a
nontrivial object: we will prove that the difference $Z^{\beta
}_t(A)-\widetilde{Z}^{\beta}_t(A)$ almost surely goes to $0$ as $t\to
\infty$ and that $\widetilde{Z}^{\beta}_t(A)$ coincides with
$M_{t}'(A)$ for $\beta$ large enough. In particular, we will prove that
$Z^{\beta}_t(A)$ converges toward a random variable $Z^{\beta}(A)$
which itself converges as \mbox{$\beta\to\infty$} to the limit of
$M_{t}'(A)$ (as $t \to\infty$). The details and proofs are gathered in
the \hyperref[app]{Appendix}.

As a direct consequence of our method of proof, we get the following
properties of $M'(dx)$:

\begin{proposition}\label{moments}
The positive random measure $M'(dx)$ possesses moments of order $q$ for
all $ q \leq0$. It does not possess moments of order $1$.
\end{proposition}

\begin{pf}
As a direct consequence of the fact that the measure
$M'$ satisfies the $\star$-equation, it possesses moments of order $q$
for all $ q\leq0$. This is a straightforward adaptation of the
corresponding theorem in \cite{Barral1}; see also \cite{Benj} for a
proof in English. Since $Z^\beta(dx)$ increases toward $M'$ as $\beta$
goes to infinity, we have $M'(dx)\geq Z^\beta(dx)$ for any $\beta$.
Since $Z^\beta_t$ is a uniformly integrable martingale, we have $
\mathbb E
[Z^\beta(A)]= \mathbb E [Z^\beta_0(A)]=\beta|A|$, we deduce that $
\mathbb E
[M'(A)]=+\infty$ for every bounded open set $A$.
\end{pf}

\section{Conjectures}\label{conj}
In this section, we present a few results we can prove about the $\star
$-equation and some conjectures related to these results.

\subsection{About the \texorpdfstring{$\star$}{star}-equation}\label{conjstar}
Consider the $\star$-equation in great generality, that is:
%
\begin{definition}[(Log-normal $\star$-scale invariance)] \label{def2}
A random Radon measure~$M$~is lognormal $\star$-scale invariant if for
all $0<\varepsilon\leq1$, $M $ obeys the cascading rule
%
\begin{equation}
\label{starc} \bigl(M (A) \bigr)_{A\in\mathcal{B}(\mathbb{R}^d)}\stackrel{\mathrm{law}} {=} \biggl(\int
_Ae^{ \omega_{\varepsilon}(r)}M^\varepsilon(dr)
\biggr)_{A\in\mathcal
{B}(\mathbb{R}^d)},
\end{equation}
where $\omega_{\varepsilon}$ is a stationary stochastically continuous
Gaussian process, and $M^{ \varepsilon}$ is a random measure
independent from $\omega_{\varepsilon}$ satisfying the scaling relation
%
\begin{equation}
\label{star1c} \bigl(M^{\varepsilon}(A) \bigr)_{A\in\mathcal{B}(\mathbb
{R}^d)}\stackrel{\mathrm{law}} {=}
\biggl(M \biggl(\frac{A}{\varepsilon}\biggr) \biggr)_{A\in\mathcal{B}(\mathbb{R}^d)}.
\end{equation}
\end{definition}
Observe that, in comparison with~(\ref{star}) and~(\ref{star1}), we do
not require the scaling factor to be $\varepsilon^d$. As stated in
(\ref{starc}) and~(\ref{star1c}), it is proved in \cite{allez} that
$ \mathbb E
[e^{ \omega_{\varepsilon}(r)}]=\varepsilon^d$ as soon as the measure
possesses a moment of order $1+\delta$ for some $\delta>0$. Roughly
speaking, it remains to investigate situations when the measure does
not possess a moment of order $1$, and we will see that the scaling
factor is then not necessarily~$\varepsilon^d$.

Inspired by the discrete multiplicative cascade case (see \cite
{durrett}), our conjecture is that all the nontrivial short ranged
solutions [i.e., there exists $R>0$ such that $M(A)$ and $M(B)$ are
independent when $d(A,B) \geq R$ where $d$ is the standard distance
between sets] to this equation belong to one of the families we will
describe below.

First we conjecture that there exists a $\alpha\in\,]0,1]$ such that
\[
\mathbb E \bigl[e^{\alpha\omega_{\varepsilon}(r)}\bigr]=\varepsilon^{ d}.
\]
Assuming this, it is proved in \cite{allez,sohier} that the Gaussian
process $\alpha\omega_{e^{-t}}$ has a covariance structure given by
(\ref{corrX}). More precisely, there exists some compactly supported
continuous covariance kernel $k$ with $k(0)=1$ and $\gamma^2\leq2d$
such that
\[
\operatorname{Cov} \bigl(\alpha\omega_{e^{-t}}(0),\alpha\omega _{e^{-t}}(x)
\bigr)=\gamma^2\int_1^{e^t}
\frac{k(ux)}{u} \,du.
\]
We can then rewrite the process $\omega$ as
\[
\omega_{e^{-t}}(x)=\frac{\gamma}{\alpha}X_t(x)-
\frac{\gamma
^2}{2\alpha
}t-\frac{d}{\alpha}t,
\]
where $(X_t)_t$ is the family of Gaussian fields introduced in
Section~\ref{setup}.
We now consider four cases, depending on the values of $\alpha$ and
$\gamma$ [cases (2), (3), (4) are conjectures]:

\begin{longlist}[(4)]
\item[(1)] If $\alpha=1$ and $\gamma^2<2d$, then the law of the
solution $M$ is the standard Gaussian multiplicative chaos $M^\gamma$
[see~(\ref{defchaos})] up to a multiplicative constant. This case has
been treated in \cite{allez}.
\item[(2)] If $\alpha=1$ and $\gamma^2=2d$, then the law of the
solution $M$ is that of the derivative martingale that we have
constructed in this paper (Theorem~\ref{mainderiv}), up to a
multiplicative constant.
\item[(3)] If $\alpha<1$ and $\gamma^2<2d$, then $M$ is an atomic
Gaussian multiplicative chaos as constructed in \cite{BJRV} up to a
multiplicative constant. More precisely, the law can be constructed as follows:
\begin{enumerate}[(a)]
\item[(a)] Sample a standard Gaussian multiplicative chaos
\[
\widebar{M}(dx)=e^{ \gamma X(x)-(\gamma^2/2) \mathbb E [X(x)^2]} \,dx.
\]
The measure $\widebar{M}$ is perfectly defined since $\gamma^2<2d$.
\item[(b)] Sample\vspace*{1pt} an independently scattered random measure $N$ whose law,
conditioned on $\widebar{M}$, is characterized by
\[
\forall q\geq0\qquad \mathbb E \bigl[e^{-qN(A)}|\widebar{M}
\bigr]=e^{-q^\alpha
\widebar{M}(A)}.
\]
Then the law of $M$ is that of $N$ up to a multiplicative constant.
\end{enumerate}
\item[(4)] If $\alpha<1$ and $\gamma^2=2d$, then $M$ is an atomic
Gaussian multiplicative chaos of a new type. More precisely, the law
can be constructed as follows:
\begin{enumerate}[(a)]
\item[(a)] Sample the derivative Gaussian multiplicative chaos
\[
M'(dx)=\bigl(\sqrt{2d} \mathbb E \bigl[X(x)^2\bigr]-
X(x)\bigr)e^{\sqrt{2d} X(x)-d
\mathbb E [X(x)^2]} \,dx.
\]
The measure $M'$ is constructed as prescribed by Theorem~\ref{mainderiv}.
\item[(b)] Sample an independently scattered random measure $N$ whose law,
conditioned on $M'$, is characterized by
\[
\forall A\in\mathcal{B}\bigl(\mathbb{R}^d\bigr),\ \forall q\geq0\qquad
\mathbb E \bigl[e^{-qN(A)}|M'\bigr]=e^{-q^\alpha M'(A)}.
\]
Then the law of $M$ is that of $N$ up to a multiplicative constant.
\end{enumerate}
\end{longlist}
Notice that the results of our paper together with \cite{allez,BJRV}
allow us to prove existence of all the random measures described above.
Therefore, it remains to complete the uniqueness part of this statement.

\begin{remark}
The $\alpha<1, \gamma^2<2d$ case above has been used in \cite
{BJRV,Duphouches,PRL} to give a mathematical understanding of the
duality in Liouville quantum gravity: this corresponds to taking
special values of the couple $(\alpha,\gamma)$. More precisely, we
choose some parameter $\bar{\gamma}^2>2d$. If the measure $M_{\bar{\gamma}}$ was well defined, it would satisfy the scaling relation
%
\begin{eqnarray}\label{stardual}
&& \bigl(M_{\bar{\gamma}}(A) \bigr)_{A\in\mathcal{B}(\mathbb{R}^d)}
\nonumber\\[-4pt]\\[-12pt]
&&\qquad \stackrel{\mathrm{law}} {=}
\biggl(\int_Ae^{ \bar{\gamma}X_{\ln(1/\varepsilon)}(r)-(\bar{\gamma}^2/2) \mathbb E [X_{\ln(1/\varepsilon)}(r)^2]}\varepsilon^d
M^{\bar{\gamma},\varepsilon}(dr) \biggr)_{A\in\mathcal{B}(\mathbb{R}^d)},\nonumber
\end{eqnarray}
where $M^{\bar{\gamma},\varepsilon}$ is a random measure independent
from $X_{\varepsilon}$ satisfying the scaling relation
%
\begin{equation}
\label{star1dual} \bigl(M^{\bar{\gamma},\varepsilon}(A) \bigr)_{A\in\mathcal
{B}(\mathbb{R}
^d)}\stackrel{\mathrm{law}} {=}
\biggl(M^{\bar{\gamma}}\biggl(\frac{A}{\varepsilon
}\biggr) \biggr)_{A\in\mathcal{B}(\mathbb{R}^d)}.
\end{equation}
Nevertheless, we know that $M^{\bar{\gamma}}$ yields a vanishing
measure. The idea is thus to use the $\star$-equation to determine what
the unique solution of this scaling relation is. Writing $\gamma=\frac
{2d}{\bar{\gamma}}<2d$ and $\alpha=\frac{2d}{\bar{\gamma}^2}$, it is
plain to see that
\[
\mathbb E \bigl[ \bigl(e^{ \bar{\gamma}X_{\ln(1/\varepsilon)}(r)-(\bar{\gamma}^2/2) \mathbb E [X_{\ln(1/\varepsilon)}(r)^2]}\varepsilon^d
\bigr)^\alpha \bigr]=\varepsilon^d.
\]
Therefore, we are in situation 3, which yields a natural candidate for
Liouville duality \cite{BJRV,Duphouches,PRL}.
\end{remark}

\subsection{Another construction of solutions to the critical \texorpdfstring{$\star$}{star}-equation}\label{secother}

Recall that the measures $M^\gamma$ for $\gamma< 2$ are obtained as
limits of~(\ref{measintro}) as $X$ varies along approximations to a
limit field. The measure constructed in Theorem~\ref{mainderiv} is
defined analogously except that one replaces~(\ref{measintro}) with
(\ref{derivintro}), which is minus the derivative of~(\ref{measintro})
at $\gamma=\sqrt{2d}$. If we could exchange the order of the
differentiation and the limit-taking, we would conclude that the
measure constructed in Theorem~\ref{mainderiv} is equal to
\[
- \frac{\partial}{\partial\gamma} \bigl[M^\gamma\bigr]_{\gamma=\sqrt{2d}} = \lim
_{\gamma\to\sqrt{2d}} \frac{1}{\sqrt{2d}-\gamma} M^\gamma.
\]

We will not fully justify this order exchange here, but we will
establish a somewhat weaker result. Namely, we show that one can at
least obtain \textit{some} solution to the $\star$-equation as a limit of
this general type. This construction is inspired by a similar
construction for discrete multiplicative cascades in \cite{durrett}.
More precisely, we have the following (proved in Section~\ref{subsecotherconsproof}):

\begin{proposition}\label{propothercons}
There exist two increasing sequence $(\lambda_n)_n$ and $(\gamma
_n)_n$, with $\gamma_n^2<2d$ and $\gamma_n^2\to2d$ as $n\to\infty$,
such that
\[
\lambda_nM^{\gamma_n}(dx)\stackrel{\mathrm{law}} {\to}
M^c(dx),
\]
where $M^c$ is a positive random measure satisfying~(\ref{star}).
\end{proposition}
The following conjecture is a consequence of the uniqueness conjecture
for the $\star$-equation exposed in Section~\ref{conjstar} above:
%
\begin{conjecture}
The construction of Proposition~\ref{propothercons} gives the same
measure as the one described in Section~\ref{secder} (up to some
multiplicative constant). Moreover, the sequence $(\lambda_n)_n$ can be
chosen as $\lambda_n=\frac{1}{\sqrt{2d}-\gamma_n}$ (in dimension~$d$).
\end{conjecture}

%
%

\subsection{Glassy phase of log-correlated Gaussian potentials}\label{glassy}
The glassy phase of log-correlated Gaussian potentials is concerned
with the behavior of
measures beyond the critical value $\gamma^2>2d$. More precisely, for
$\gamma^2>2d$, consider the
measure
\[
M^{\gamma}_t(dx)=e^{\gamma X_t(x)-(\gamma^2/2) \mathbb E
[X_t(x)^2]} \,dx.
\]
The limiting measure, as $t\to\infty$, vanishes as proved in \cite
{cfKah}. Therefore, it is natural to look for a suitable family of
normalizing factors to make this measure converge. With the arguments
used in Section~\ref{cascade} to compare with the results obtained
in \cite{BRV,madaule}, we can rigorously prove:
%
\begin{proposition}
The renormalized family
\[
\bigl( t^{(3\gamma)/(2\sqrt{2d})}e^{t ((\gamma/\sqrt{2})-\sqrt{d} )^2}M^\gamma_t(dx)
\bigr) _{t \geq0}
\]
is tight. Furthermore, every converging subsequence is nontrivial.
\end{proposition}
The above proposition can be obtained using the results in \cite
{BRV,madaule} and Section~\ref{cascade} (tightness statement). The main
result in \cite{bramson} about the behavior of the maximum of the
discrete GFF implies that every converging subsequence is nontrivial.

We now formulate a conjecture about the limiting law of this family.
Assuming that the above renormalized family converges in law (so we
strengthen tightness into convergence), it turns out that the limit
$M^\gamma$ of this renormalized family necessarily satisfies the
following $\star$-equation:
\[
M^\gamma(dx)=e^{ \gamma X_{\ln(1/\varepsilon)}(x)-\sqrt{(d/2)} \gamma \mathbb E [X_{\ln(1/\varepsilon)}(x)^2]} \varepsilon ^{\sqrt{(d/2)} \gamma}
\widebar{M}^{\gamma} \biggl(\frac
{dx}{\varepsilon} \biggr),
\]
where $\widebar{M}^{\gamma}$ is a random measure with the same law as
$M^\gamma$ and independent of the process $(X_t(x))_{x\in\mathbb{R}^d}$.
Setting $\alpha=\frac{\sqrt{2d}}{\gamma}\in\,]0,1[$, this equation
can be
rewritten as
\[
M^\gamma(dx)=e^{(\sqrt{2d}/\alpha) X_{\ln(1/\varepsilon)}(x)-(d/\alpha) \mathbb E [X_{\ln(1/\varepsilon)}(x)^2]} \varepsilon^{d/\alpha}
\widebar{M}^{\gamma} \biggl(\frac
{dx}{\varepsilon
} \biggr).
\]
Therefore, assuming that the conjectures about uniqueness of the $\star
$-equation are true, we have the following:

\begin{conjecture}
%
\begin{equation}
\label{renormsurcrit} t^{(3\gamma)/(2\sqrt{2d})}e^{t ((\gamma/\sqrt{2})-\sqrt
{d} )^2}M^\gamma_t(dx)
\stackrel{\mathrm{law}} {\to} c_\gamma N_\alpha(dx)\qquad\mbox{as }t\to
\infty,
\end{equation}
where $c_\gamma$ is a positive constant depending on $\gamma$ and the
law of $N_\alpha$ is given, conditioned on the derivative martingale
$M'$, by an independently scattered random measure the law of which is
characterized by
\[
\forall A\in\mathcal{B}\bigl(\mathbb{R}^d\bigr),\ \forall q\geq0\qquad
\mathbb E \bigl[e^{-qN_\alpha
(A)}|M'\bigr]=e^{-q^\alpha M'(A)}.
\]
\end{conjecture}

In particular, physicists are interested in the behavior of the Gibbs
measure associated to $M^{\gamma}_t(dx)$ on a ball $B$. It is the
measure renormalized by its total mass,
\[
G^\gamma_t(dx)=\frac{M^{\gamma}_t(dx)}{M^{\gamma}_t(B)}.
\]
From~(\ref{renormsurcrit}), we deduce
%
\begin{equation}
G^\gamma_t(dx)\stackrel{\mathrm{law}} {\to} \frac{N_\alpha(dx)}{N_\alpha
(B)}\qquad
\mbox{as }t\to\infty.
\end{equation}
The size reordered atoms of the latter object form a Poisson--Dirichlet
process as conjectured by physicists \cite{CarDou} and proved
rigorously in \cite{arguin}. Nevertheless, we point out that this
conjecture is more powerful than the Poisson--Dirichlet result since it
also makes precise the spatial localization of the atoms. We stress
that this result has been proved in the case of branching random walks
\cite{BRV}, built on the work of Madaule \cite{madaule}.


\subsection{About the maximum of the log-correlated Gaussian random variables}\label{sec4.4}
It~is proved in \cite{bramson} (in fact $d=2$ in \cite{bramson} but
this is general) that the family
\[
\biggl(\sup_{x \in[0,1]^d} X_t(x)- \sqrt{2d} t +
\frac{3}{2\sqrt {2d}} \ln t \biggr)_{t \geq0}
\]
is tight. One can thus conjecture by analogy with the branching random
walk case~\cite{aidekon}:
%
\begin{conjecture}
\[
\sup_{x \in[0,1]^d} X_t(x)- \sqrt{2d} t +
\frac{3}{2\sqrt{2d}} \ln t \to G_d\qquad\mbox{in law as }t\to\infty,
\]
where the distribution of $G_d$ is given in terms of the distribution
of the limit $M'([0,1]^d)$ of the derivative martingale. More
precisely, there exists some constant $c>0$ such that
%
\begin{equation}
\label{eqonmax} \mathbb E \bigl[e^{-qG_d}\bigr]= \frac{1}{c^q} \Gamma
\biggl(1+\frac{q}{\sqrt{2d}}\biggr) \mathbb E \bigl[ \bigl( M'
\bigl([0,1]^d\bigr) \bigr) ^{-q/\sqrt{2d}} \bigr].
\end{equation}
\end{conjecture}
%
Here we give a heuristic derivation of identity~(\ref{eqonmax}) using
the conjectures of the above subsections. By performing an inversion of limits:
($\gamma\leftrightarrow t$ and conjecturing $\frac{\ln c_\gamma
}{\gamma}\to\ln c$ as $\gamma\to\infty$),
\begin{eqnarray*}
&& \mathbb E \bigl[e^{-qG_d}\bigr]
\\
&&\qquad =\lim_{\gamma\to+ \infty}\lim
_{t \to+ \infty} \mathbb E \bigl[\exp \bigl[-q{\gamma}^{-1}
\ln \bigl[t^{(3\gamma)/(2\sqrt{2 d})}e^{t  ((\gamma/\sqrt{2})-\sqrt{d} )^2} M^\gamma _t
\bigl([0,1]^d\bigr) \bigr] \bigr] \bigr]
\\
&&\qquad =\lim_{\gamma\to+ \infty} \mathbb E \bigl[ \bigl(c_\gamma N_{\alpha=\sqrt{2d}/\gamma}\bigl([0,1]^d\bigr) \bigr)^{-q/\gamma}  \bigr]
\\
&&\qquad =\frac{1}{c^q} \Gamma\biggl(1+\frac{q}{\sqrt{2d}}\biggr) \mathbb E
\bigl[{ \bigl( M'\bigl([0,1]^d\bigr) \bigr)
^{-q/\sqrt{2d}} } \bigr],
\end{eqnarray*}
%
where, for $x>0$, $\Gamma(x)=\int_0^{\infty} t^{x-1} e^{-t} \,dt $ is the
standard Gamma function. Therefore, $G_d$ can be viewed as a modified
Gumbel law. Otherwise stated, we conjecture
\begin{eqnarray*}
&& \lim_{t\to\infty} \mathbb P \biggl( \sup_{x \in[0,1]^d}
X_t(x)- \sqrt{2d} t +\frac{3}{2\sqrt{2d}} \ln t\leq u \biggr)
\\
&&\qquad = \mathbb
E \bigl[\exp \bigl[-c^{\sqrt{2d}}e^{-
\sqrt{2d} u}M'
\bigl([0,1]^d\bigr)\bigr] \bigr].
\end{eqnarray*}
We point out that we recover in a heuristic and alternative way the
result proved rigorously in \cite{aidekon} for branching random walks.

\begin{appendix}\label{app}
\section{Proofs}\label{sec5}
\subsection{Proofs of results from Section~\texorpdfstring{\protect\ref{secder}}{3}}

We follow the notation of Section~\ref{secder}. We first investigate
the convergence of $(Z^\beta_t(A))_{t\geq0}$:

\begin{proposition}
The process $(Z^\beta_t(A))_{t\geq0}$ is a continuous positive
\mbox{$\mathcal{F}_t$-}mar\-tingale and thus converges almost surely toward a
positive random variable denoted by $Z^\beta(A)$. 
\end{proposition}

\begin{pf}
Proving that $(Z^\beta_t(A))_{t\geq0}$ is a
martingale boils down to proving, for each $x\in A$, that
\begin{eqnarray*}
&& \mathbb E \bigl[\bigl(\sqrt{2d}t-X_{t}(x)+\beta\bigr)
\one_{\{\tau^\beta> t\}}e^{\sqrt
{2d}X_{t}(x)-d \mathbb E [X_{t}(x)^2]}|\mathcal{F}_s\bigr]
\\
&&\qquad =\bigl(\sqrt
{2d}s-X_{s}(x)+\beta \bigr)\one_{\{\tau^\beta> s\}}e^{\sqrt{2d}X_{s}(x)-d \mathbb E [X_{s}(x)^2]}.
\end{eqnarray*}
Let us first stress that, for each $x\in A $, the process
$(X_t(x))_{t\geq0}$ is a Brownian motion. Furthermore, we can use the
(weak) Markov property of the Brownian motion to get
\begin{eqnarray*}
&& \mathbb E \bigl[\bigl(\sqrt{2d}t-X_{t}(x)+\beta\bigr)
\one_{\{\tau^\beta> t\}
}e^{\sqrt
{2d}X_{t}(x)-d \mathbb E [X_{t}(x)^2]} |\mathcal{F}_s\bigr]
\\
&&\qquad =\one_{\{\tau^\beta> s\}}e^{\sqrt{2d}X_{s}(x)-d \mathbb E
[X_{s}(x)^2]}F \bigl(\sqrt{2d}s-X_{s}(x)+
\beta \bigr),
\end{eqnarray*}
where
\begin{eqnarray*}
F(y)&=& \mathbb E \bigl[\bigl(\sqrt{2d}(t-s)-X_{t-s}(x)+y\bigr)
\\
&&{}\hspace*{9pt}\times \one_{\{\tau(X_\cdot
(x)-\sqrt {2d}\cdot-y)> t-s\}} e^{\sqrt{2d}X_{t-s}(x)-d E[X_{t-s}(x)^2]} \bigr]
\end{eqnarray*}
and, for a stochastic process $Y$, $\tau(Y)$ is defined by
\[
\tau(Y)=\inf\{u>0;Y_u>0\}.
\]
Using the Girsanov transform yields
\[
F(y)= \mathbb E \bigl[\bigl(-X_{t-s}(x)+y\bigr)\one_{\{\tau(X_\cdot(x)-y)> t-s\}}
\bigr].
\]
Hence we get
\begin{eqnarray*}
F(y) &=& \mathbb E \bigl[\bigl(-X_{t-s}(x)+y\bigr)\one_{\{\tau(X_\cdot(x)-y)> t-s\}}
\bigr]
\\
&=& \mathbb E \bigl[\bigl(-X_{(t-s) \wedge\tau(X_\cdot(x)-y)}(x)+y\bigr) \bigr]=y
\end{eqnarray*}
by the optional stopping theorem. This completes the proof.
\end{pf}

\begin{proposition}\label{proplog}
Assume that $A$ is a bounded open set. Then the martingale $(Z^\beta
_t(A))_{t\geq0}$ is regular.
\end{proposition}

\begin{pf}
Without loss of generality, we may assume
$k(u)=0$ for $|u|>1$ since $k$ has a compact support (so we just assume
that the smallest ball centered at $0$ containing the support of $k$
has radius $1$ instead of $R$ for some $R>0$). We may also assume that
$A\subset B(0,1/2)$: indeed, any bigger bounded set can be recovered
with finitely many balls with radius less than $\frac{1}{2}$. Finally,
we will also assume that $x\cdot\nabla k(x)\leq0$. This condition need
not be true over the whole $\mathbb{R}^d$. Nevertheless, it must be
valid in a
neighborhood of $0$ [and even $x\cdot\nabla k(x)< 0 $ if $x\neq0$] in
order not to contradict the fact that $k$ is positive definite and nonconstant. Therefore, even if it means considering a smaller set $A$, we
may (and will) assume that this condition holds.

Write for $x\in\mathbb{R}^d$
\[
f_t^\beta(x)=\bigl(\sqrt{2d}t-X_{t}(x)+\beta
\bigr)\one_{\{\tau^\beta>t\}} e^{\sqrt
{2d}X_{t}(x)-d \mathbb E [X_{t}(x)^2]}.
\]
Define then the analog of the \textit{rooted random measure} in \cite{cfDuSh} (also called the ``Peyri\`ere probability measure'' in this
context \cite{cfKah}),
\[
\Theta^\beta_t=\frac{1}{|A|\beta}f_t^\beta(x)\,dx
\,d \mathbb P.
\]
It is a probability measure on $\mathcal{B}(A)\otimes\mathcal{F}_t $.
We denote by $\Theta^\beta_t(\cdot|\mathcal{G})$ the conditional
expectation of $\Theta^\beta_t$ given some sub-$\sigma$-algebra
$\mathcal{G}$ of $\mathcal{B}(A)\otimes\mathcal{F}_t $. If $y$ is a
$\mathcal{B}(A)\otimes\mathcal{F}_t $-measurable random variable on
$A\times\Omega$, we denote by $\Theta^\beta_t(\cdot|y)$ the
conditional expectation of $\Theta^\beta_t$ given the $\sigma$-algebra
generated by $y$.

We first observe that
\[
\Theta^\beta_t(\cdot|x) =\frac{1}{\beta}f_t^\beta(x)\,d \mathbb P.
\]
Therefore, under $\Theta^\beta_t(\cdot|x)$, the process
$(X_{s}(x)-\sqrt {2d}s-\beta)_{s\leq t}$ has the law of $(-\beta_s)_{s\leq t}$ where
$(\beta_s)_{s\leq t}$ is a $3d$-Bessel process starting from $\beta$.
Let us now recall the following result (see \cite{Motoo}):

\begin{theorem}\label{thmotoo}
Let $X$ be a $3d$-Bessel process on $\mathbb{R}_+$ started from $
\beta\geq0$
with respect to the law $ \mathbb P _\beta$.
\begin{longlist}[(2)]
\item[(1)] Suppose that $\phi\uparrow\infty$ such that $\int_1^{\infty
}\frac{\phi(t)^3}{t}e^{-(1/2)\phi(t)^2} \,dt<+\infty$. Then
\[
\mathbb P _\beta \bigl(X_t>\sqrt{t}\phi(t)\mbox{ i.o.
as } t\uparrow+\infty \bigr)=0.
\]
\item[(2)] Suppose that $\psi\downarrow0$ such that $\int_1^{\infty
}\frac{\psi(t)}{t}  \,dt<+\infty$. Then
\[
\mathbb P _\beta \bigl(X_t<\sqrt{t}\psi(t)\mbox{ i.o.
as } t\uparrow+\infty \bigr)=0.
\]
\end{longlist}
\end{theorem}
In view of the above theorem, we can choose $R$ large enough such that
for all $x$ the set
\[
B_t(x)= \biggl\{\forall s \in[0,t]; \frac{\sqrt{s}}{R(\ln(2+s))^2} \leq\beta+
\sqrt{2d}s-X_{s}(x) \leq R\bigl(1+\sqrt{s\ln(1+s)}\bigr) \biggr\}
\]
has a probability arbitrarily close to $1$, say $1-\varepsilon$, for
all $t\dvtx  \Theta^\beta_t(B_t(x) |x) \geq1-\varepsilon$.

We can now prove the uniform integrability of $(Z^\beta_t(A))_t$, that is,
\[
\lim_{\delta\to\infty}\limsup_{t\to\infty} \mathbb E
\bigl[Z^\beta _t(A)\one_{\{Z^\beta
_t(A)>\delta\}}\bigr]=0.
\]
Observe that
\[
\mathbb E \bigl[Z^\beta_t(A)\one_{\{Z^\beta_t(A)>\delta\}}\bigr]=
\beta |A|\Theta^\beta _t\bigl(Z^\beta_t(A)>
\delta\bigr).
\]
Therefore, it suffices to prove that
\[
\lim_{\delta\to\infty}\limsup_{t\to\infty}
\Theta^\beta_t \bigl(Z_t^\beta (A)>
\delta\bigr)=0.
\]
We have
\begin{eqnarray*}
&& \Theta^\beta_t \bigl(Z_t^\beta(A)>
\delta\bigr)
\\
&&\qquad = \frac{1}{|A|}\int_A\Theta^\beta_t
\bigl(Z_t^\beta(A)>\delta|x\bigr) \,dx
\\
&&\qquad = \frac{1}{|A|}\int_A\Theta^\beta_t
\bigl( \Theta^\beta_t \bigl(Z_t^\beta(A)>
\delta|x,\bigl(X_s(x)\bigr)_{s\leq t} \bigr) |x \bigr) \,dx
\\
&&\qquad \leq\varepsilon+ \frac{1}{|A|}\int_A
\Theta^\beta_t \bigl( \Theta ^\beta_t
\bigl(Z_t^\beta(A)>\delta|x,\bigl(X_s(x)
\bigr)_{s\leq t},B_t(x) \bigr) |x \bigr) \,dx
\\
&&\qquad \leq\varepsilon+ \frac{1}{|A|}\int_A
\Theta^\beta_t \biggl(\Theta ^\beta _t
\biggl(Z_t^\beta\bigl(B\bigl(x,e^{-t}\bigr)
\bigr)>\frac{\delta}{2} \bigg|x,\bigl(X_s(x)\bigr)_{s\leq
t},B_t(x) \biggr)\bigg|x \biggr) \,dx
\\
&&\quad\qquad{} + \frac{1}{|A|}\int_A\Theta^\beta_t
\biggl(\Theta^\beta_t \biggl(Z_t^\beta
\bigl(B\bigl(x,e^{-t}\bigr)^c\bigr)>\frac{\delta}{2} \bigg|x,
\bigl(X_s(x)\bigr)_{s\leq t},B_t(x) \biggr)\bigg|x
\biggr) \,dx
\\
&&\qquad \stackrel{\mathrm{def}} {=}\varepsilon+\Pi_1+\Pi_2.
\end{eqnarray*}

We are now going to estimate $\Pi_1,\Pi_2$. Observe that the two
quantities roughly reduce to expressions like ($K$ is a ball or its
complementary)
\[
\Theta^\beta_t \biggl(H \biggl(\int_{K}f_t^\beta(w)
\,dw \biggr) \bigg|x,\bigl(X_s(x)\bigr)_{s\leq t},B_t(x)
\biggr)
\]
for $H$ a nonnegative function (here an indicator function). To carry
out our computations, we thus have to compute the law of the process
$(X_s(w))_{s\leq t}$ knowing that of the process $(X_s(x))_{s\leq t}$.
To that purpose, we will use the following lemma whose proof is left to
the reader since it follows from a standard (though not quite direct)
computation of covariances for Gaussian processes:
%
\begin{lemma}\label{decomp}
For $w\neq x$ and all $s_0$, the law of the process $(X_s(w))_{s\leq s_0}$
can be decomposed as
\[
X_s(w)=P^{x,w}_s+Z^{x,w}_s,
\]
where:
\begin{longlist}[--]
\item[--] $P^{x,w}_s=-\int_0^sg_{x,w}(u)X_u(x) \,du+ K'_s(x-w) X_s(x)$ is
measurable with respect to the $\sigma$-algebra generated by
$(X_s(x))_{s\leq s_0}$ and $g_{x,w}(u)= K''_u(x-w) $;

\item[--]the process $(Z^{x,w}_s)_{0\leq s\leq s_0}$ is a centered Gaussian
process independent of $(X_s(x))_{0\leq s\leq s_0}$ with covariance kernel
\[
q_{x,w}\bigl(s,s'\bigr)\stackrel{\mathrm{def}} {=} \mathbb E
\bigl[Z^{x,w}_{s'}Z^{x,w}_s\bigr]= s
\wedge s' -\int_0^{s\wedge s'}
\bigl(K'_u(x-w)\bigr)^2 \,du.
\]
\end{longlist}%
\end{lemma}
The above decomposition lemma roughly implies the following: the two
processes $(X_s(w))_{s\geq0}$ and $(X_s(x))_{s\geq0}$ are the same
until $s_0=\ln\frac{1}{|x-w|}$ and then the two processes
$(X_s(w)-X_{s_0}(w))_{s\geq s_0}$ and $(X_s(x)-X_{s_0}(x))_{s\geq s_0}$
are independent.

We first estimate $\Pi_2$ with the above lemma. It is enough to
estimate properly the quantity
%
\begin{equation}
\label{pi2tilde} \widetilde{\Pi}_2=\Theta^\beta_t
\biggl(Z_t^\beta \bigl(B\bigl(x,e^{-t}
\bigr)^c\bigr)>\frac
{\delta}{2} \bigg|x,\bigl(X_s(x)
\bigr)_{s\leq t},B_t(x) \biggr).
\end{equation}
Notice that
%
\begin{equation}
\label{pi2tilde2} \widetilde{\Pi}_2\leq\frac{2}{\delta}\int
_{B(x,e^{-t})^c}\Theta ^\beta _t
\bigl(f^\beta_t(w) |x,\bigl(X_s(x)
\bigr)_{s\leq t},B_t(x) \bigr) \,dw.
\end{equation}
For each $w\in B(x,e^{-t})^c $, that is, such that $|w-x|>e^{-t}$, let
us define $s_0=\ln\frac{1}{|x-w|}$. Notice\vspace*{-2pt} that $s_0$ is
the time at which the evolution of $(X_s(w)-X_{s_0}(w))_{s_0\leq s\leq
t}$ becomes independent of the process $(X_s(x))_{0\leq s\leq t}$.
Under $\Theta^\beta_t$, the process $(X_s(w))_{s_0\leq s\leq t}$ can be
rewritten as
\[
X_s(w)=X_{s_0}(w)+W_{s-s_0},
\]
where $W$ is a standard Brownian motion independent of the
processes\break
$(X_s(x))_{0\leq s\leq t}$ and $(X_s(w))_{0\leq s\leq s_0}$. This can
be checked by a straightforward computation of covariance.
Therefore, we get
\begin{eqnarray*}
&&\Theta^\beta_t\bigl(f_t^\beta(w)|x,
\bigl(X_s(x)\bigr)_{s\leq t}\bigr)
\\[-1pt]
&&\qquad =\frac{1}{\beta} \mathbb E \bigl[\bigl(\sqrt{2d}t-X_t(w)+\beta
\bigr)
\\[-1pt]
&&\hspace*{53pt}\times{} \one_{\{\sup_{[0,t]}X_u(w)-\sqrt{2d}u\leq\beta\}} e^{\sqrt
{2d}X_t(w)-dt}|x,\bigl(X_s(x)
\bigr)_{s\leq t} \bigr]
\\[-1pt]
&&\qquad =\frac{1}{\beta} \mathbb E \bigl[\bigl(\sqrt{2d}s_0+\sqrt
{2d}(t-s_0)-X_{s_0}(w)-W_{t-s_0}+\beta\bigr)
\\[-1pt]
&&\hspace*{53pt}\times{} \one_{\{\sup
_{[0,s_0]}X_u(w)-\sqrt{2d}u\leq\beta\}}
\\[-1pt]
&&\hspace*{53pt}\times{} \one_{\{\sup_{[s_0,t]}X_{s_0}(w)+\sqrt{2d}s_0+W_{u-s_0}-\sqrt
{2d}(u-s_0)\leq\beta\}}
\\[-1pt]
&&\hspace*{53pt}\times{}  e^{\sqrt{2d}X_{s_0}(w)-ds_0} e^{\sqrt
{2d}W_{t-s_0}-d(t-s_0)} |x,
\bigl(X_s(x)\bigr)_{s\leq t} \bigr]
\\[-1pt]
&&\qquad =\frac{1}{\beta} \mathbb E \bigl[\bigl(\sqrt{2d}s_0-X_{s_0}(w)+
\beta \bigr)
\\[-1pt]
&&\hspace*{53pt}\times{} \one_{\{\sup
_{[0,s_0]}X_u(w)-\sqrt{2d}u\leq\beta\}} e^{\sqrt{2d}X_{s_0}(w)-ds_0} |x,\bigl(X_s(x)
\bigr)_{s\leq t} \bigr]
\end{eqnarray*}
by the stopping time theorem. From Lemma~\ref{decomp}, we deduce
%
\begin{eqnarray}\label{estpi2}
\hspace*{-2pt}&& \Theta^\beta_t \bigl(f_t^\beta(w)|x,
\bigl(X_s(x)\bigr)_{s\leq t}\bigr)\nonumber
\\
\hspace*{-2pt}&&\qquad =\frac{1}{\beta} \mathbb E \bigl[\bigl(\sqrt {2d}s_0-P^{x,w}_{s_0}-Z^{x,w}_{s_0}+
\beta \bigr)\nonumber
\\
\hspace*{-2pt}&&\hspace*{51pt}{}\times \one_{\{\sup_{[0,s_0]}P^{x,w}_{u}+Z^{x,w}_{u}-\sqrt{2d}u\leq\beta\}}
e^{\sqrt{2d}(P^{x,w}_{s_0}+Z^{x,w}_{s_0})-ds_0} |x,\bigl(X_s(x)
\bigr)_{s\leq t} \bigr]\hspace*{-15pt}\nonumber
\\
\hspace*{-2pt}&&\qquad \leq\frac{1}{\beta} \mathbb E \bigl[ \bigl(\bigl(\sqrt {2d}s_0-P^{x,w}_{s_0}-Z^{x,w}_{s_0}+
\beta\bigr)^2+1 \bigr)
\\
\hspace*{-2pt}&&\hspace*{51pt}{}\times  e^{\sqrt
{2d}(P^{x,w}_{s_0}+Z^{x,w}_{s_0})-ds_0} |x,\bigl(X_s(x)
\bigr)_{s\leq t} \bigr]\nonumber
\\
\hspace*{-2pt}&&\qquad =\frac{1}{\beta} \bigl(\bigl(\sqrt {2d}\bigl(s_0-q_{x,w}(s_0,s_0)
\bigr)-P^{x,w}_{s_0}+\beta\bigr)^2+q_{x,w}(s_0,s_0)
\bigr)\nonumber
\\
\hspace*{-2pt}&&\hspace*{30pt}{} \times  e^{\sqrt{2d}P^{x,w}_{s_0} -d(s_0-q_{x,w}(s_0,s_0))}.\nonumber
\end{eqnarray}
We make two observations. First, we point out that the quantity
$q_{x,w}(s_0,s_0)$ is bounded by a constant only depending on $k$ since
\begin{eqnarray*}
q_{x,w}(s_0,s_0)&=&s_0 -\int
_0^{s_0}\bigl(K'_u(x-w)
\bigr)^2 \,du
\\
&=& \int_0^{s_0} \bigl[1- \bigl(k
\bigl(e^u(x-w)\bigr) \bigr)^2 \bigr] \,du
\\
&= & \int_{|x-w|}^{1} \biggl(1- k \biggl(y
\frac{x-w}{|x-w|} \biggr) ^2 \biggr) \frac
{1}{y} \,dy
\\
&\leq&C,
\end{eqnarray*}
where $C$ can be defined as $\sup_{z\in B(0,1)}\frac{1-k(z)^2}{|z|} $.
So the quantity $q_{x,w}(s_0,s_0)$ will not play a part in the forthcoming
computations.

Second, we want to express the random variable $P^{x,w}_{s_0}$ as a
function of the Bessel process $(X_u(x)-\sqrt{2d}u-\beta)_u$ in order
to use the fact that we can control the paths of this latter process
[we will condition by the event $B_t(x)$]. Therefore, we set
%
\begin{eqnarray}\label{defY}
Y_{s_0}^{x,w}&=&-\int_0^{s_0}g_{x,w}(u)
\bigl(X_u(x)-\sqrt{2d}u-\beta\bigr) \,du\nonumber
\\
&=&-\int_0^{s_0}g_{x,w}(u)X_u(x)
\,du-\sqrt{2d}K_{s_0}(x-w)
\nonumber\\[-8pt]\\[-8pt]
&&{} +\beta \bigl(k\bigl(e^{s_0}(x-w)
\bigr)-k(x-w) \bigr)\nonumber
\\
&=&P^{x,w}_{s_0}-\sqrt{2d}K_{s_0}(x-w)+\beta \bigl(k
\bigl(e^{s_0}(x-w)\bigr)-k(x-w) \bigr).
\nonumber
\end{eqnarray}
Therefore, we can write
\[
Y_{s_0}^{x,w}=P^{x,w}_{s_0}-
\sqrt{2d}s_0+\theta_{x,w}(s_0)
\]
for some function $\theta_{x,w}$ that is bounded independently of
$x,w,t$ since $k$ is bounded over $\mathbb{R}^d$. Plugging these
estimates into
(\ref{estpi2}), we obtain
%
\begin{eqnarray}\label{majbess}
&&\Theta^\beta_t\bigl(f_t^\beta(w)|x,
\bigl(X_s(x)\bigr)_{s\leq t}\bigr)
\nonumber
\\
&&\qquad =\frac{1}{\beta} \bigl(\bigl(\theta _{x,w}(s_0)-Y^{x,w}_{s_0}
\bigr)^2+q_{x,w}(s_0,s_0) \bigr)
\nonumber\\[-8pt]\\[-8pt]
&&\hspace*{30pt}{}\times e^{\sqrt{2d}Y^{x,w}_{s_0}+ds_0+dq_{x,w}(s_0,s_0)-\sqrt{2d}\theta _{x,w}(s_0)} \nonumber
\\
&&\qquad \leq\frac{C}{\beta} \bigl(\bigl(Y^{x,w}_{s_0}
\bigr)^2+1\bigr) e^{\sqrt
{2d}Y^{x,w}_{s_0}+ds_0}\nonumber
\end{eqnarray}
for some constant $C$ that does not depend on $x,w,t$. Now we plug the
exact expression of $g_{x,w}$,
\[
g_{x,w}(u)=\sum_{i=1}^d(x-w)_ie^u
\partial_i k\bigl(e^u(x-w)\bigr)
\]
into definition~(\ref{defY}) of $Y^{x,w}_{s_0}$,
\begin{eqnarray*}
Y^{x,w}_{s_0}&=&\int_0^{\ln(1/|x-w|)}
\sum_{i=1}^d(x-w)_ie^u
\partial _i k\bigl(e^u(x-w)\bigr) \bigl(\sqrt{2d}u+
\beta-X_u(x) \bigr) \,du
\\
&=&\int_{|x-w|}^{1}y\frac{x-w }{|x-w|}\cdot\nabla k
\biggl(y\frac{x-w
}{|x-w|} \biggr)
\\
&&\hspace*{27pt}{}\times \biggl(\sqrt{2d}\ln\frac{y}{|x-w|}+
\beta-X_{\ln(y/|x-w|)}(x) \biggr) \,dy.
\end{eqnarray*}
Moreover the constraint for the Bessel process, valid on $B_t(x)$,
%
\begin{eqnarray}\label{strip}
\frac{\sqrt{u}}{R(\ln(2+u))^2} &\leq&\beta-
X_u(x)+\sqrt{2d}u \leq R\bigl(1+\sqrt{u\ln(1+u)}\bigr)
\nonumber\\[-8pt]\\[-12pt]
\eqntext{\forall u \in[0,t]}
\end{eqnarray}
implies that [here we use the relation $x\cdot\nabla k(x)\leq0$]
%
\begin{eqnarray}
\qquad Y^{x,w}_{s_0} & \geq& R\int_{|x-w|}^{1}y
\frac{x-w }{|x-w|}\cdot\nabla k \biggl(y\frac{x-w }{|x-w|} \biggr)
\nonumber\\[-8pt]\\[-8pt]
&&\hspace*{36pt}{}\times  \biggl(1+
\sqrt{\ln\frac{y}{|x-w|}\ln \biggl(1+\ln\frac{y}{|x-w|} \biggr)} \biggr) \,dy,\nonumber
\\
Y^{x,w}_{s_0} &\leq& R\int_{|x-w|}^{1}y
\frac{x-w }{|x-w|}\cdot\nabla k \biggl(y\frac{x-w }{|x-w|} \biggr)
\frac{\sqrt{\ln(y/|x-w|)}}{\ln(2+\ln(y/|x-w|))^2} \,dy.
\end{eqnarray}
Using rough estimates yields
%
\begin{eqnarray}
&& -C_R \biggl(1+\sqrt{\ln\frac{1}{|x-w|}\ln \biggl(1+\ln
\frac{1}{|x-w|} \biggr)} \biggr) \,du
\nonumber\\[-8pt]\\[-8pt]
&&\qquad \leq Y^{x,w}_{s_0} \leq-C_R \frac{\sqrt{\ln(1/|x-w|)}}{\ln (2+\ln(1/|x-w|) )^2}\nonumber
\end{eqnarray}
for some constant $C_R$ depending on $R$ and on the function $x\mapsto
x\cdot\nabla k(x)$. Plugging these estimates into~(\ref{majbess})
yields (the constant $C$ may change, depending on the value of $C_R$)
%
\begin{equation}
\label{jenaichie} \Theta^\beta_t \bigl(f_t^\beta(w)|x,
\bigl(X_s(x)\bigr)_{s\leq t},B_t(x)\bigr)\leq
\frac
{e^C}{\beta|x-w|^d}G\biggl(\ln\frac{1}{|x-w|}\biggr),
\end{equation}
where
\[
G(y)= \bigl(1+\sqrt{y\ln (1+y )} \bigr)^2 e^{-\sqrt{2d}C(\sqrt{y}/\ln (2+y )^2)}.
\]
Finally, by gathering estimates~(\ref{pi2tilde}),~(\ref{pi2tilde2}) and
(\ref{jenaichie}) and then making successive changes of variables, we
obtain ($V_d$ stands for the area of the unit sphere of $\mathbb{R}^d$)
\begin{eqnarray*}
\Pi_2&=&\frac{1}{|A|}\int_A
\Theta^\beta_t \biggl(\Theta^\beta _t
\biggl(Z_t^\beta\bigl(B\bigl(x,e^{-t}
\bigr)^c\bigr)>\frac{\delta}{2} \bigg|x,\bigl(X_s(x)
\bigr)_{s\leq t},B_t(x) \biggr)\bigg|x \biggr) \,dx
\\
&=&\frac{1}{|A|}\int_A\Theta^\beta_t
(\widetilde{\Pi}_2|x ) \,dx
\\
&\leq&\frac{2}{|A|\delta}\int_A\int_{B(x,e^{-t})^c}
\frac
{e^C}{\beta
|x-w|^d}G\biggl(\ln\frac{1}{|x-w|}\biggr)\,dx\,dw
\\
&\leq&\frac{2V_d}{ \delta}\int_{e^{-t}}^{1}
\frac{e^C}{\beta
r^d}G\biggl(\ln \frac{1}{r}\biggr)r^{d-1} \,dr
\\
&\leq&\frac{2V_de^C}{ \delta\beta}\int_{0}^{t} G(u) \,du.
\end{eqnarray*}
Since $G$ is integrable, this quantity is obviously bounded by a
quantity that goes to~$0$ when $\delta$ becomes large uniformly with
respect to $t$. This concludes estimating~$\Pi_2$.

We now estimate $\Pi_1$. Once again, it is enough to estimate the quantity
%
\begin{equation}
\label{pi1tilde} \widetilde{\Pi}_1=\Theta^\beta_t
\biggl(Z_t^\beta\bigl(B\bigl(x,e^{-t}\bigr)
\bigr)>\frac{\delta}{2} \bigg|x,\bigl(X_s(x)\bigr)_{s\leq t},B_t(x)
\biggr),
\end{equation}
which is less than
%
\begin{equation}
\label{pi1tilde2} \widetilde{\Pi}_1\leq\frac{2}{\delta}\int
_{B(x,e^{-t}) }\Theta ^\beta _t
\bigl(f^\beta_t(w) |x,\bigl(X_s(x)
\bigr)_{s\leq t},B_t(x) \bigr) \,dw.
\end{equation}
This time, for $|w-x|<e^{-t}$, there is no need to ``cut'' the process
$(X_s(w))_{s\leq t}$ at level $s_0=\ln\frac{1}{|x-w|}$. We can
directly use Lemma~\ref{decomp} to get
\begin{eqnarray*}
\hspace*{-5pt}&&\Theta^\beta_t\bigl(f_t^\beta(w)|x,
\bigl(X_s(x)\bigr)_{s\leq t},B\bigr)
\\
\hspace*{-5pt}&&\qquad =\frac{1}{\beta} \mathbb E \bigl[\bigl(\sqrt {2d}t-P^{x,w}_{t}-Z^{x,w}_{t}+
\beta\bigr)
\\
\hspace*{-5pt}&&\hspace*{51pt}{}\times \one _{\{\sup_{[0,t]}P^{x,w}_{u}+Z^{x,w}_{u}-\sqrt{2d}u\leq\beta\}}
\\
\hspace*{-5pt}&&\hspace*{51pt}{}\times e^{\sqrt
{2d}(P^{x,w}_{t}+Z^{x,w}_{t})-dt} |x,\bigl(X_s(x)
\bigr)_{s\leq t},B_t(x) \bigr]
\\
\hspace*{-5pt}&&\qquad \leq\frac{1}{\beta} \mathbb E \bigl[ \bigl(\bigl(\sqrt {2d}t-P^{x,w}_{t}-Z^{x,w}_{t}+
\beta\bigr)^2+1 \bigr)
\\
\hspace*{-5pt}&&\hspace*{51pt}{}\times  e^{\sqrt
{2d}(P^{x,w}_{t}+Z^{x,w}_{t})-dt} |x,\bigl(X_s(x)
\bigr)_{s\leq t},B_t(x) \bigr]
\\
\hspace*{-5pt}&&\qquad =\frac{1}{\beta} \bigl(\bigl(\sqrt{2d}\bigl(t-q_{x,w}(t,t)
\bigr)-P^{x,w}_{t}+\beta \bigr)^2+q_{x,w}(t,t)
\bigr) e^{\sqrt{2d}P^{x,w}_{t} -d(t-q_{x,w}(t,t))}.
\end{eqnarray*}
Once again, the quantity $q_{x,w}(t,t)$ is bounded by a constant only
depending on $k$ (not on $t$).
Second, for $s\leq t$, we define the process
\begin{eqnarray*}
Y_{s}^{x,w}&=&-\int_0^{s }g_{x,w}(u)
\bigl(X_u(x)-\sqrt{2d}u-\beta\bigr) \,du
\\
&&{} +K'_s(x-w) \bigl(X_s(x)-\sqrt{2d}s-\beta\bigr),
\end{eqnarray*}
which turns out to be equal to
\begin{eqnarray*}
Y_{s}^{x,w}&=&P^{x,w}_{s}-\sqrt{2d}s+
\theta_{x,w}(s)
\end{eqnarray*}
for some function $\theta_{x,w}$ that is bounded independently of
$x,w,s$. We deduce
%
\begin{eqnarray}\label{majbess1}
\qquad &&\Theta^\beta_t\bigl(f_t^\beta(w)|x,
\bigl(X_s(x)\bigr)_{s\leq t}\bigr)
\nonumber
\\
&&\qquad =\frac{1}{\beta} \bigl(\bigl(\theta _{x,w}(t)-Y^{x,w}_{t}
\bigr)^2+q_{x,w}(t,t) \bigr) e^{\sqrt{2d}Y^{x,w}_{t}+dt+dq_{x,w}(t,t)-\sqrt{2d}\theta
_{x,w}(t)}
\\
&&\qquad \leq\frac{C}{\beta} \bigl(\bigl(Y^{x,w}_{t}\bigr)^2+1\bigr)  e^{\sqrt
{2d}Y^{x,w}_{t}+dt}\nonumber
\end{eqnarray}
for some constant $C$ that does not depend on $x,w,t$. Once again on
$B_t(x)$, the Bessel process evolves in the strip~(\ref{strip}),
implying that the process $Y^{x,w}$ is bound to live in the strip (we
stick to the previous notations)
%
\begin{equation}
-C_R \bigl(1+\sqrt{t\ln (1+t )} \bigr) \,du\leq Y^{x,w}_{t}
\leq-C_R \frac{\sqrt{t}}{\ln (2+ t )^2}
\end{equation}
for some constant $C_R$. Plugging these estimates into~(\ref{majbess1})
yields (the constant $C$ may change, depending on the value of $C_R$)
%
\begin{equation}
\label{jenaichie1} \Theta^\beta_t \bigl(f_t^\beta(w)|x,
\bigl(X_s(x)\bigr)_{s\leq t},B_t(x)\bigr)\leq
\frac
{e^C}{\beta}G(t)e^{dt},
\end{equation}
where the function $G$ is still defined by
\[
G(t)= \bigl(1+\sqrt{t\ln (1+t )} \bigr)^2 e^{-\sqrt{2d}C
(\sqrt{t}/\ln(2+t )^2)}.
\]
Notice that this estimate differs from that obtained for $\widetilde
{\Pi
}_2$ because of the $e^{dt}$ factor. It will be absorbed by the volume
of the ball $B(x,e^{-t})$ that we will integrate over. Finally, by
using~(\ref{jenaichie1}), we obtain
\begin{eqnarray*}
\Pi_1&=&\frac{1}{|A|}\int_A
\Theta^\beta_t \biggl(\Theta^\beta _t
\biggl(Z_t^\beta\bigl(B\bigl(x,e^{-t}\bigr)
\bigr)>\frac{\delta}{2} \bigg|x,\bigl(X_s(x)\bigr)_{s\leq t},B_t(x)
\biggr)\bigg|x \biggr) \,dx
\\
&=&\frac{1}{|A|}\int_A\Theta^\beta_t
(\widetilde{\Pi}_2|x ) \,dx
\\
&\leq&\frac{2}{|A|\delta}\int_A\int_{B(x,e^{-t}) }
\frac{e^C}{\beta
}G(t)e^{dt}\,dx\,dw
\\
&\leq&\frac{2 }{ \delta} \frac{e^C}{\beta}G(t).
\end{eqnarray*}
Since $G$ is bounded, this quantity is obviously bounded by a quantity
that goes to $0$ when $\delta$ becomes large uniformly with respect to
$t$. This concludes estimating $\Pi_1$. The proof is complete.
\end{pf}

We are now in position to prove the following:
%
\begin{theorem}
For each bounded open set $A\subset\mathbb{R}^d$, the martingale\break
$(M'_t(A))_{t\geq0}$ converges almost surely toward a positive random
variable denoted by $M'(A)$, such that $M'(A)>0$ almost surely.
Consequently, almost surely, the (locally signed) random measures
$(M'_t(dx))_{t\geq0}$ converge weakly as $t\to\infty$ toward a
positive random measure $M'(dx)$, which has full support and is
atomless. Furthermore, the measure $M'$ is a solution to the $\star
$-equation~(\ref{star}) with $\gamma=\sqrt{2d}$.
\end{theorem}

\begin{pf}
We first observe that the martingale $(Z^\beta
_t(A))_{t\geq0}$ possesses almost surely the same limit as the process
$(\widetilde{Z}^\beta_t(A))_{t\geq0}$
because
%
\begin{equation}
\label{diffz} \qquad\bigl|Z^\beta_t(A)-\widetilde{Z}^\beta_t(A)\bigr|=
\beta\int_A\one_{\{\tau
^\beta>
t\}} e^{\sqrt{2d}X_{t}(x)-dE[X_{t}(x)^2]} \,dx\leq
\beta M^{\sqrt{2d}}_t(A)
\end{equation}
and the last quantity converges almost surely toward $0$ since
$M^{\sqrt
{2d}}_t(dx)$ almost surely converges toward $0$ as $t$ goes to $\infty
$; see Proposition~\ref{propstand} below. Using Proposition \ref{propstand}, we have almost surely,
\[
\sup_{t\in\mathbb{R}_+}\max_{x\in A}X_{t}(x)-
\sqrt{2d}t<+\infty,
\]
which obviously implies
\[
\forall t\qquad M'_t(A)=\widetilde{Z}^\beta_t(A)
\]
for $\beta$ (random) large enough.

Since the family of random measures $(Z^\beta_t(dx))_{t\geq0}$ are
nonnegative, and $(Z^\beta_t(A))_{t\geq0}$ almost surely converges for
every bounded open set $A$, it is plain to deduce that, almost surely,
the random measures $(Z^\beta_t(dx))_{t\geq0}$ and $(\widetilde
{Z}^\beta_t(dx))_{t\geq0}$ weakly converge toward a random measure
$Z^\beta(dx)$. Then, almost surely, the family $(M'_t(dx))_{t\geq0}$
weakly converges toward the positive random measure defined by the
increasing limit $M'(dx):= \lim_{\beta\to\infty} Z^\beta(dx)$.
Indeed, consider $L>0$. We want to show that $(M'_t(dx))_{t\geq0}$
converges weakly on $[-L,L]^d$. If $\varepsilon>0$, we can find a
$\beta
>0$ such that
%
\begin{equation}
\label{eqevent} E_\beta(L):= \sup_{t\in\mathbb{R}_+}\max
_{x\in
[-L,L]^d}X_{t}(x)-\sqrt{2d}t \leq\beta
\end{equation}
has probability greater or equal to $1-\varepsilon$. On the event
$E_\beta(L)$, we have for all $\beta' \geq\beta$ the following equality:
\[
M'_t(A)= Z^\beta_t(A)-\beta
M^{\sqrt{2d}}_t(A),\qquad t \geq0, A \subset[-L,L]^d.
\]
Hence, on the event $E_\beta$, the signed measure $M'_t(dx)$ converges
weakly on $[-L,L]^d$ toward $M'(dx)=Z^\beta(dx)$.

Let us prove that the support of $M'$ is $\mathbb{R}^d$. We first
write the
relation, for $s<t$,
%
\begin{eqnarray}
\label{passlim} Z^\beta_t(dx)&=&\bigl(\sqrt{2d}s-X_{s}(x)+
\beta\bigr)\one_{\{\tau^\beta>t\}} e^{\sqrt{2d}X_{t}(x)-d \mathbb E [X_{t}(x)^2]} \,dx\nonumber
\\
&&{} +\bigl(\sqrt {2d}(t-s)-X_{t}(x)+X_{s}(x)+\beta\bigr)
\\
&&\hspace*{10pt}{}\times  \one_{\{\tau^\beta>t\}}e^{\sqrt{2d}X_{t}(x)-d \mathbb E [X_{t}(x)^2]} \,dx.
\nonumber
\end{eqnarray}
By using the same arguments as throughout this section, we pass to the
limit in this relation as $t\to\infty$ and then $\beta\to\infty$
to get
%
\begin{eqnarray}
\label{equiv} M' (dx)&=& e^{\sqrt{2d}X_{s}(x)-d \mathbb E [X_{s}(x)^2]} M^{\prime,s}(dx),
\end{eqnarray}
where $M^{\prime,s}$ is defined as
\[
M^{\prime,s}(dx)=\lim_{\beta\to\infty}\lim_{t \to\infty}Z^{\beta,s}_t
(dx)
\]
and $Z_t^{\beta,s}(dx)$ is almost surely defined as the weak limit of
\begin{eqnarray*}
Z^{\beta,s}_t (A) &=& \int_A \bigl(\sqrt
{2d}(t-s)-X_{t}(x)+X_{s}(x)+\beta \bigr)
\\
&&\hspace*{9pt}{}\times  \one_{\{\tau^\beta_s>t\}} e^{\sqrt{2d}(X_{t}(x)-X_{s}(x))-d (
\mathbb E [X_{t}(x)^2]- \mathbb E [X_{s}(x)^2] )} \,dx,
\end{eqnarray*}
where
\[
\tau^\beta_s=\inf\bigl\{u>0;X_{u+s}(x)-X_{s}(x)-
\sqrt{2d}u>\beta\bigr\}.
\]
Let us stress that we have used the fact that the measure
\[
\bigl(\sqrt{2d}s-X_{s}(x)+\beta\bigr)\one_{\{\tau^\beta>t\}}
e^{\sqrt
{2d}X_{t}(x)-d \mathbb E [X_{t}(x)^2]} \,dx
\]
goes to $0$ [it is absolutely continuous w.r.t. to $M^{\sqrt
{2d}}_t(dx)$] when passing to the limit in~(\ref{passlim}). Therefore,
$M'$ is a solution to the $\star$-equation~(\ref{star}). From (\ref
{equiv}), it is plain to deduce that the event $\{M'(A)=0\}$ ($A$ open
nonempty set) belongs to the asymptotic sigma-algebra generated by the
field $\{(X_t(x))_x;t\geq0\}$. Therefore, it has probability $0$ or $1$
by the $0-1$ law of Kolmogorov. Since we have already proved that it is
not $0$, this proves that $ \mathbb P (M'(A)=0)=0$ for any nonempty
open set $A$.

Finally, we prove that the measure is atomless. The proof is based on
the computations made during the proof of Proposition~\ref{proplog}.
We will explain how to optimize these computations to obtain the
atomless property. Of course, we could have done that directly in the
proof of Proposition~\ref{proplog}, but we feel that it is more
pedagogical to separate the arguments. Let us roughly explain how we
will proceed. Clearly, it is sufficient to prove that the positive
random measure
\[
Z^\beta(dx)= \lim_{t \to\infty}Z^{\beta}_t
(dx)
\]
does not possess atoms. Indeed, on the event $E_\beta(L)$ defined by
(\ref{eqevent}), the measure $M'(dx)$ coincides with $Z^\beta(dx)$ on
$[-L,L]^d$.

To that purpose, by stationarity, it is enough to prove that (see \cite{daley}, Corollary~9.3, Chapter~VI)
\[
\forall\delta>0\qquad\lim_nn^d \mathbb P
\bigl(Z^\beta(I_n)>\delta \bigr)=0,
\]
where $I_n$ is the cube $[0,\frac{1}{n}]^d$.
From now on, we stick to the notations of Proposition~\ref{proplog}.
We have to prove that
\[
\forall\delta>0\qquad\lim_n \limsup_t
\Theta_t^\beta \bigl(Z^\beta _t(I_n)>
\delta \bigr)=0.
\]

Therefore, let $\delta>0$ and $\varepsilon>0$ be two fixed positive
real numbers. We choose $R$ and the associated event $B$ of probability
$1-\varepsilon$ as in Proposition~\ref{proplog}. We have
\[
\limsup_t \Theta_t^\beta
\bigl(Z^\beta_t(I_n)>\delta \bigr) \leq
\varepsilon+ \limsup_t \Pi_1+ \limsup
_t \Pi_2.
\]
First note that $\limsup_t   \Pi_1=0$; we also have the following
bound for $\limsup_t   \Pi_2$:
\[
\limsup_t \Pi_2 \leq\frac{2V_de^C}{ \delta\beta}\int
_{n \ln2
}^{\infty} G(u) \,du,
\]
which goes to $0$ as $n$ goes to $\infty$. In conclusion, we get
\[
\lim_n \limsup_t
\Theta_t^\beta \bigl(Z^\beta_t(I_n)>
\delta \bigr) \leq \varepsilon,
\]
which is the desired result.
\end{pf}

\subsection{Proof of result from Section~\texorpdfstring{\protect\ref{conj}}{4}} \label{subsecotherconsproof}

Here, we prove Proposition~\ref{propothercons}. For notational
simplicity, we further assume that the dimension $d$ is equal to $1$
and that $k(u)=0$ for all $|u|>1$. Generalization to all other
situations is straightforward.

Let $C$ be the interval $[0,1]$. Let us denote by $\phi(\cdot,\gamma)$
the Laplace transform of~$M^\gamma(C)$
\[
\phi(\lambda,\gamma)= \mathbb E \bigl[e^{-\lambda M^\gamma(C)}\bigr].
\]
Since $ \mathbb P (M^\gamma(C)>0)=1$ the range of the mapping $\lambda
\in\mathbb{R}
_+\mapsto\phi(\lambda,\gamma)$ is the whole interval $]0,1]$.
Choose a
strictly increasing sequence $(\gamma_n)_n$ converging toward $\sqrt {2}$. Choose a sequence $(\lambda_n)_n$ such that
%
\begin{equation}
\label{tight} \phi(\lambda_n,\gamma_n)=
\tfrac{1}{2}.
\end{equation}
Let us denote by $M^c(C)$ a random variable taking values in
$[0,+\infty
]$ such that $\lambda_nM^{\gamma_n}(C)\to M^c(C)$ vaguely as $n\to
\infty$ (eventually up to a subsequence). Let us define the function
\[
\varphi(\theta)= \mathbb E \bigl[e^{-\theta M^c(C)},M^c(C)<\infty
\bigr]
\]
for $\theta>0$ and $\varphi(0)=1$. Then $\phi(\theta\lambda
_n,\gamma
_n)\to\varphi(\theta)$ for all $\theta$ so that, in particular,
$\varphi(1)=\frac{1}{2}$. Let us choose $\varepsilon$ small enough in
order to have $\ln\frac{1}{\varepsilon}$ even integer larger than $4$.
Because of~(\ref{star}), we have
\begin{eqnarray*}
\phi(\theta\lambda_n,\gamma_n) &=& \mathbb E \biggl[
\exp \biggl[{-\theta\lambda _n\int_Ce^{ \gamma_n X_{\ln(1/\varepsilon)}(r)-(\gamma
_n^2/2) \mathbb E [X_{\ln(1/\varepsilon)}(r)^2]}M^{\gamma
_n,\varepsilon
}(dr)}
\biggr] \biggr].
\end{eqnarray*}
Let us denote by $C_k$ the interval $[\frac{k}{\ln(1/\varepsilon)},\frac{k+1}{\ln(1/\varepsilon)}]$ for $k\in A_\varepsilon
\stackrel{\mathrm{def}}{=}\{0,\ldots, \ln\frac{1}{\varepsilon}-1\}$. By the
Cauchy--Schwarz inequality and stationarity, we have
\begin{eqnarray*}
\phi(\theta\lambda_n,\gamma_n) &\leq& \mathbb E \biggl[
\exp \biggl[-2\theta\lambda _n\mathop{\sum
_{k \in A_\varepsilon}}_{\mathrm{even}}\int_{C_k}e^{
\gamma_n X_{\ln(1/\varepsilon)}(r)-(\gamma_n^2/2)\mathbb E
[X_{\ln(1/\varepsilon)}(r)^2]}M^{\gamma_n,\varepsilon
}(dr)
\biggr] \biggr].
\end{eqnarray*}
By the Kahane convexity inequality and because the mapping $x\mapsto
e^{-sx}$ is convex for any $s\in\mathbb{R}$, we deduce
\begin{eqnarray*}
\phi(\theta\lambda_n,\gamma_n) &\leq& \mathbb E \biggl[
\exp \biggl[-2\theta\lambda _n\mathop{\sum
_{k \in A_\varepsilon}}_{\mathrm{even}}\int_{C_k}e^{
\sqrt{2}X_{\ln(1/\varepsilon)}(0)- \mathbb E [X_{\ln(1/\varepsilon)}(0)^2]}M^{\gamma_n,\varepsilon}(dr)
\biggr] \biggr]
\\
&=& \mathbb E \biggl[\exp \biggl[{-2\theta\lambda_ne^{ \sqrt{2} X_{\ln
(1/\varepsilon)}(0)- \mathbb E [X_{\ln(1/\varepsilon)}(0)^2]}
\mathop{\sum_{k \in A_\varepsilon}}_{\mathrm{even}}
M^{\gamma
_n,\varepsilon}(C_k)} \biggr] \biggr].
\end{eqnarray*}
Because\vspace*{1pt} the sets $C_k$ are separated by a distance of at least $ \frac
{1}{\ln(1/{\varepsilon})}$, the random variables $\mathop{(M^{\gamma
_n,\varepsilon}(C_k))_{k \in A_\varepsilon\ \mathrm{even}}}$
are i.i.d. with common law $\varepsilon M^{\gamma_n }(C )$ because of~(\ref{star}). We deduce
\begin{eqnarray*}
\phi(\theta\lambda_n,\gamma_n) &\leq& \mathbb E \bigl[
\phi \bigl(2\theta\lambda _n\varepsilon e^{ \sqrt{2} X_{\ln(1/\varepsilon)}(0)- \mathbb
E [X_{\ln(1/\varepsilon)}(0)^2]},
\gamma_n \bigr)^{(1/2)\ln(1/\varepsilon)} \bigr].
\end{eqnarray*}
By taking the limit as $n\to\infty$, we deduce
\begin{eqnarray*}
\varphi(\theta)&\leq& \mathbb E \bigl[\varphi \bigl(2\theta \varepsilon
e^{ \sqrt
{2} X_{\ln(1/\varepsilon)}(0)- \mathbb E [X_{\ln(1/\varepsilon)}(0)^2]} \bigr)^{(1/2)\ln(1/\varepsilon)} \bigr].
\end{eqnarray*}
By letting $\theta$ go to $0$, we deduce
\[
\varphi(0_+)\leq\varphi(0_+)^{(1/2)\ln(1/\varepsilon)}.
\]
Because $\frac{1}{2}\ln\frac{1}{\varepsilon}\geq2$, we are left with
two options: either $\varphi(0_+)=0$ or \mbox{$\varphi(0_+)\geq1$}. But
$\varphi(0_+)\leq1$ because $e^{-\theta x}\leq1$ for all $x\geq0$.
Furthermore $\varphi(0_+)\geq\varphi(1)=\frac{1}{2}$. Therefore,
$\varphi(0_+)=1$ and $M^c(C)<+\infty$ almost surely. $M^c(C)$ is not
trivial because $\varphi(1)=\frac{1}{2}$. We have proved that the
sequence $(\lambda_nM^{\gamma_n}(C))_n$ is tight and that the limit of
every converging subsequence is nontrivial.

Of course, we can carry out the same job for every smaller dyadic
interval. But the normalizing sequence may depend on the size of the
interval. Let us prove that it does not. To this purpose, it is enough
to establish that
\[
\frac{1}{2}\leq\liminf_{n} \mathbb E
\bigl[e^{-\lambda_n M^{\gamma
_n}(C_k)} \bigr]\leq \limsup_{n} \mathbb E
\bigl[e^{-\lambda_n M^{\gamma_n}(C_k)} \bigr]<1
\]
for every dyadic interval $C_k$ of size $2^{-k}$. The left-hand side is
obvious because $M^{\gamma_n}(C_k)\leq M^{\gamma_n}(C)$. By using
(\ref{star}) with $\varepsilon=2^{-k}$ and the Kahane convexity inequality,
we deduce
\begin{eqnarray*}
&& \limsup_{n} \mathbb E \bigl[\exp \bigl[{-\lambda_n M^{\gamma
_n}(C_k)} \bigr]\bigr]
\\
&&\qquad \leq \limsup_{n}
\mathbb E \bigl[\exp \bigl[{-\lambda_n M^{\gamma
_n}(C)2^{-k}e^{\sqrt
{2}X_{k\ln2}(0)- \mathbb E [X_{k\ln2}(0)^2]}}
\bigr] \bigr]
\\
&&\qquad = \mathbb E \bigl[\varphi \bigl(2^{-k}e^{\sqrt{2}X_{k\ln2}(0)-
\mathbb E [X_{k\ln
2}(0)^2]} \bigr)
\bigr].
\end{eqnarray*}
The last quantity is strictly less than $1$. Indeed, if not, then
\[
\varphi \bigl(2^{-k}e^{\sqrt{2d}X_{k\ln2}(0)-((2d)/2) \mathbb E
[X_{k\ln
2}(0)^2]} \bigr)=1
\]
almost surely, that is, $\varphi(\theta)=1$ for all
$\theta$, hence a contradiction.

To sum up, the sequence $(\lambda_n M^{\gamma_n}(C))_n$ is tight for
all dyadic intervals. By the Tychonoff theorem and the Caratheodory
extension theorem, we can extract a subsequence and find a random
measure $M^c(dx)$ such that $(\lambda_n M^{\gamma_n}(C_1),\ldots,\lambda
_n M^{\gamma_n}(C_p))_n$ converges in law toward $(
M^{c}(C_1),\ldots,\break
M^{c}(C_p))_n$ for all dyadic intervals $C_1,\ldots,C_p$. Finally, by
multiplying both sides of~(\ref{star}) by $\lambda_n$ and passing to
the limit as $n\to\infty$, we deduce
%
\begin{equation}
\label{stargamc} \quad \bigl(M^c(A) \bigr)_{A\in\mathcal{B}(\mathbb{R})}\stackrel{\mathrm{law}} {=}
\biggl(\int_Ae^{ \sqrt{2 } X_{\ln(1/\varepsilon)}(r)- \mathbb E [X_{\ln(1/\varepsilon)}(r)^2]}M^{c,\varepsilon}(dr)
\biggr)_{A\in\mathcal
{B}(\mathbb{R})},
\end{equation}
where
%
\begin{equation}
\label{star1gamc} \bigl(M^{c,\varepsilon}(A) \bigr)_{A\in\mathcal{B}(\mathbb
{R})}\stackrel {\mathrm{law}}
{=}\varepsilon \biggl(M^c\biggl(\frac{A}{\varepsilon}\biggr)
\biggr)_{A\in
\mathcal
{B}(\mathbb{R})}.
\end{equation}
%

\section{Auxiliary results}\label{sec6}
We first state the classical ``Kahane's convexity inequalities''
(originally written in \cite{cfKah}; see also \cite{allez} for a proof):

\begin{lemma}\label{lemcvx}
Let $F,G\dvtx \mathbb{R}_+\to\mathbb{R}$ be two functions such that $F$
is convex, $G$ is
concave and
\[
\forall x\in\mathbb{R}_+\qquad\bigl|F(x)\bigr|+\bigl|G(x)\bigr|\leq M\bigl(1+|x|^\beta
\bigr)
\]
for some positive constants $M,\beta$, and $\sigma$ be a Radon measure
on the Borelian subsets of $\mathbb{R}^d$. Given a bounded Borelian
set $A$,
let $(X_r)_{r\in A},(Y_r)_{r\in A}$ be two continuous centered Gaussian
processes with continuous covariance kernels $k_X$~and~$k_Y$ such that
\[
\forall u,v\in A\qquad k_X(u,v)\leq k_Y(u,v).
\]
Then
\begin{eqnarray*}
\mathbb E \biggl[F \biggl(\int_A e^{X_r-(1/2) \mathbb E [X_r^2]}
\sigma(dr) \biggr) \biggr]&\leq& \mathbb E \biggl[F \biggl(\int_Ae^{Y_r-(1/2) \mathbb E
[Y_r^2] }\sigma(dr) \biggr) \biggr],
\\
\mathbb E \biggl[G \biggl(\int_A e^{X_r-(1/2) \mathbb E [X_r^2]}
\sigma(dr) \biggr) \biggr]&\geq& \mathbb E \biggl[G \biggl(\int_Ae^{Y_r-(1/2) \mathbb E
[Y_r^2] }\sigma(dr) \biggr) \biggr].
\end{eqnarray*}
If we further assume
\[
\forall u\in A\qquad k_X(u,u)= k_Y(u,u),
\]
then we recover Slepian's comparison lemma: for each increasing
function $F\dvtx \mathbb{R}_+\to\mathbb{R}$:
\[
\mathbb E \Bigl[F \Bigl(\sup_{x\in A}Y_x \Bigr)
\Bigr]\leq \mathbb E \Bigl[F \Bigl(\sup_{x\in
A}X_x
\Bigr) \Bigr].
\]
\end{lemma}

\subsection{Chaos associated to cascades}\label{cascade}

We use Kahane convexity inequalities (see Proposition~\ref{lemcvx}) to
compare the small moments of the Gaussian multiplicative chaos with
those of a dyadic lognormal Mandelbrot's multiplicative cascade. Let us
briefly recall the construction of lognormal Mandelbrot's
multiplicative cascades. We consider the $2^d$-adic tree
\[
T=\bigl(\{1,2\}^d\bigr)^{\mathbb{N}^*}.
\]
For $t\in T$, we denote by $\pi_k(t)$ ($k\in\mathbb{N}^*$) the $k$th
component
of $t$. We equip $T$ with the ultrametric distance
\[
\forall s,t\in T\qquad\mathbf{d}(t,s)=2^{-dn}\qquad\mbox{where } n=
\sup\bigl\{N\in\mathbb{N}; \forall k\leq N, \pi_k(t)=
\pi_k(s)\bigr\}
\]
with the convention that $n=0$ if the set $\{N\in\mathbb{N}; \forall
k\leq N,
\pi_k(t)=\pi_k(s)\}$ is empty.
Let us define
\[
\forall s,t \in T\qquad p_n(t,s)=\cases{ u, &\quad if $
\mathbf{d}(t,s)\leq2^{-nd}$,
\vspace*{2pt}\cr
0, &\quad if $\mathbf{d}(t,s)>
2^{-nd}$.}
\]
The kernel $p_n$ is therefore constant over each of the $2^{dn}$
cylinders defined by the prescription of the first $n$ coordinates [in
what follows, we will denote by $I_n(t)$ that cylinder containing $t$].
For each $n$, we denote by $(Y_n(t))_{t\in T}$ a centered Gaussian
process indexed by $T$ with covariance kernel $p_n$. We assume that the
processes~$(Y_n)_n$ are independent.
We set
%
\begin{equation}
\label{kernelq} \forall s,t \in T\qquad q_n(t,s)=\sum
_{k=1}^np_k(t,s).
\end{equation}
Notice that
%
\begin{equation}
\label{kernelexpl} \forall s,t \in T\qquad q_n(t,s)= \frac{u}{d \ln2}
\ln\frac
{1}{\mathbf{d}(t,s)\vee2^{-dn}}
\end{equation}
and
\[
q_n(t,s)\to\frac{u}{d \ln2}\ln\frac{1}{\mathbf{d}(t,s)} \qquad \mbox{as }n\to\infty.
\]
We define the centered Gaussian process
\[
\forall t\in T\qquad \widebar{X}_n(t)=\sum
_{k=1}^nY_k(t)
\]
with covariance kernel $q_n$. Let us denote by $\sigma$ the uniform
measure on $T$, that is $\sigma(I_n(t))=2^{-dn}$. We set
\[
\widebar{M}^u_n=\int_Te^{\widebar{X}_n(t)-(1/2) \mathbb E
[\widebar{X}_n(t)^2]}
\sigma(dt).
\]
This corresponds to the lognormal multiplicative cascades framework.
The martingale $(\widebar{M}^u_n)_n$ converges toward a nontrivial
limit if and only if $u<2d \ln2$. The boundary case corresponds to
$u=2d \ln2$. It is proved in \cite{cfKahPer} that, for $u=2d \ln2$,
$\lim_n\widebar{M}^u_n(dx)=0$ almost surely.

It turns out that the $2^d$-adic tree can be naturally embedded in the
unit cube of $\mathbb{R}^d$ by iteratively dividing a cube into $2^d$ cubes
with equal size length. 
Notice that the uniform measure on the tree is then sent to the
Lebesgue measure by this embedding. We also stress that the dyadic
distance on the cube $[0,1]^d$ is greater than the Euclidean distance
on that cube
\[
\forall s,t\in[0,1]^d\qquad|t-s| \leq\sqrt{d} \mathbf{d}(t,s)^{1/d}.
\]
This allows many one-sided comparison results between lognormal
cascades and Gaussian multiplicative chaos.

So, taking $u= 2d \ln2$ in the kernel $q_n$ of~(\ref{kernelexpl}), we
claim for all $s',s\in[0,1]^d,\ \forall n\in\mathbb{N}$,
%
\begin{equation}
\label{compcasc} q_n\bigl(s,s'\bigr)-C\leq2d
K_{n\ln2}\bigl(s-s'\bigr)
\end{equation}
for some constant $C>0$ that does not depend on $n$ (only on $k$).

We are now in position to prove the following:
%
\begin{proposition}\label{propstand}
For $\gamma^2=2d$, the standard construction yields a vanishing
limiting measure
%
\begin{equation}
\label{eqngamma2limitzero} \lim_{t \to\infty} M_t^{\sqrt{2d}} =0
\qquad\mbox{almost surely}.
\end{equation}
Furthermore, for all $a \in[0,\frac{1}{2}[$ and any bounded open set
$A$, almost surely,
%
\begin{equation}
\label{max} \sup_{t \geq0} \biggl(\sup_{x\in A}X_{t}(x)-
\sqrt{2d}t+\frac{a}{\sqrt{2d}} \ln(t+1)\biggr) < \infty.
\end{equation}
\end{proposition}

\begin{pf}
We consider $\widebar{X}_n$ with covariance given by (\ref
{kernelexpl}) for $u=\ln2$; by a slight abuse of notation, we consider
that $\widebar{X}_n$ is defined on the unit cube by the natural embedding.

The family $(M^\gamma_t)$ is a positive martingale. Therefore, it
converges almost surely. We just have to prove that the limit is zero.
We will apply Kahane's concentration inequalities (Lemma \ref
{lemcvx}). Let us denote by $Z$ a standard Gaussian random variable
independent of the process $(X_t(x))_{t,x}$. From~(\ref{compcasc}),
the covariance kernel of the centered Gaussian process $\widebar{X}_n$
is less than that of the Gaussian process $\sqrt{C}Z+X_{n\ln2}$. By
applying Lemma~\ref{lemcvx} to some bounded concave function
$F\dvtx \mathbb{R}_+\to
\mathbb{R}$ and $n\in\mathbb{N}$, we obtain (we stick to the
notations introduced just above)
%
\begin{equation}
\label{compccv} \mathbb E \bigl[F \bigl( e^{\sqrt{C}Z-(1/2)C}M^{\sqrt
{2d}}_{n\ln
2}
\bigl([0,1]^d\bigr) \bigr) \bigr]\leq \mathbb E \biggl[F \biggl( \int
_{T} e^{\sqrt
{2d}
\widebar{X}_n(t)-d  \mathbb E [\widebar{X}_n(t)^2]} \,dt \biggr) \biggr].\hspace*{-35pt}
\end{equation}
Now we further assume that $F$ is increasing. Because of the dominated
convergence theorem, the right-hand side goes to $F(0)$ as $n\to\infty
$. So does the left-hand side. This shows that $M^{\sqrt{2d}}_{n\ln
2}([0,1]^d)$ goes to $0$ in probability as $n\to\infty$. Since we
already know that the martingale $M^{\sqrt{2d}}_{t}([0,1]^d)$ converges
almost surely as $t\to\infty$, this completes the proof of the first
statement.

For the second statement, we fix $a \in[0,\frac{1}{2}[$ and we
consider the case $d=1$ with $k(x)=(1-|x|)_{+}$ for simplicity\vadjust{\goodbreak} (this is
no restriction since every $C^1$ kernel $k$ with $k(0)=1$ is greater or
equal to some $(1-\frac{|x|}{L})_{+}$ for $L>0$). In this case, one can
represent the variables $X_s(x)$ as integrals of truncated cones with
respect to a Gaussian measure; see Section~\ref{cone} below for a
quick reminder or \cite{Bar,bacry} for details. Note that a similar
cone construction can be performed in $\mathbb{R}^d \times\mathbb
{R}_+$, and hence the
proof can be generalized to all dimensions. The cone representation
ensures that we have the following decomposition (see Section~\ref{cone}):

\begin{lemma}\label{lemcone}
We fix $n$ and cut $[0,1]$ into $2^n$ intervals. We have the following
decomposition for $X_{s \ln2}(x)$ for all $s \in[n,n+1]$ and $x \in
I_{i,n}:=[\frac{i}{2^n},\frac{i+1}{2^n}[$:
\[
X_{s \ln2}(x)=X_{i,n}+Y_{s}^{i,n}(x)
\]
with the following properties:
\begin{itemize}
\item
There exists a constant $C>0$ (independent of $n$) such that
\begin{eqnarray*}
\mathbb E [X_{i,n} X_{j,n}] &=& n \ln2-\biggl(1-
\frac{1}{2^n}\biggr)\qquad\mbox{if } i=j,
\\[-1pt]
\mathbb E [X_{i,n} X_{j,n}] &\geq& E\biggl[
\widebar{X}_{n}\biggl(\frac{i}{2^n}\biggr) \widebar{X}_{n}
\biggl(\frac{j}{2^n}\biggr)\biggr]-C\qquad\mbox{if } i \neq j.
\end{eqnarray*}
\item
For all $i$, the process $(Y_{s}^{i,n}(x))_{s \in[n,n+1], x \in
I_{i,n}}$ is continuous and independent of~$X_{i,n}$.
\item
For all $i,j$, $s,s' \in[n,n+1]$ and $x \in I_{i,n}$, $x' \in I_{j,n}$:
\[
\mathbb E \bigl[ Y_{s}^{i,n}(x) Y_{s'}^{j,n}
\bigl(x'\bigr) \bigr] \geq0.
\]
\item
For all $i,j$, $s \in[n,n+1]$ and $x \in I_{i,n}$:
\[
\mathbb E \bigl[ Y_{s}^{i,n}(x) X_{j,n} \bigr]
\geq0.
\]
\end{itemize}
\end{lemma}

We\vspace*{1pt} introduce a standard Gaussian variable $Z$ independent from the
process $(X_{s \ln2}(x))_x$ and a standard Gaussian\vspace*{1pt} i.i.d. sequence
$(\widebar{Z}_{i})_{0 \leq i \leq2^n-1}$. We also introduce a
sequence of independent processes $(\widebar{Y}_{s}^{i,n}(x))_{s \in
[n,n+1], x \in I_{i,n}}$ independent from $\widebar{X}_{n}$
and such that for all $i$ the process $(\widebar{Y}_{s}^{i,n}(x))_{s
\in[n,n+1], x \in I_{i,n}}$ has same law as $(Y_{s}^{i,n}(t))_{s \in
[n,n+1], x \in I_{i,n}}$.
By Lemma~\ref{lemcvx}, we have the following for all $y$:
\begin{eqnarray*}
&& \mathbb P \biggl( \sup_{0 \leq i \leq2^n-1}\,\sup_{s \in[n,n+1]}\,\sup_{x \in
I_{i,n}} \biggl( X_{i,n}+\sqrt{1-
\frac{1}{2^n}+C}Z+Y_{s}^{i,n}(x) -\sqrt {2} n \ln2
\biggr) \geq y \biggr)
\\[-1pt]
&&\qquad  \leq \mathbb P \biggl( \sup_{0 \leq i \leq2^n-1}\,\sup_{s \in[n,n+1]}\,\sup_{x \in
I_{i,n}} \biggl(\widebar{X}_{n}\biggl(
\frac{i}{2^n}\biggr)+\sqrt{C}\, \widebar{Z}_i
\\[-3pt]
&&\hspace*{155pt}{} + \widebar{Y}_{s}^{i,n}(x) - \sqrt{2} n \ln2\biggr) \geq y \biggr).
\end{eqnarray*}
Indeed, we have the following if $i=j$, $x,x' \in I_{i,n}$ and $s,s'
\in[n,n+1]$:
\begin{eqnarray*}
&& \mathbb E \biggl[ \biggl(X_{i,n}+\sqrt{1-\frac{1}{2^n}+C}Z+Y_{s}^{i,n}(x)
\biggr) \biggl(X_{i,n}+\sqrt {1-\frac{1}{2^n}+C}Z+Y_{s'}^{i,n}
\bigl(x'\bigr)\biggr)\biggr]
\\[-1pt]
&&\qquad  = n \ln2+C + \mathbb E \bigl[ Y_{s}^{i,n}(x)
Y_{s'}^{i,n}\bigl(x'\bigr) \bigr]
\\[-1pt]
&&\qquad =\mathbb E \biggl[\biggl(\widebar{X}_{n}\biggl(
\frac{i}{2^n}\biggr)+\sqrt{C} \widebar{Z}_i+
\widebar{Y}_{s}^{i,n}(x)\biggr) \biggl(
\widebar{X}_{n}\biggl(\frac
{i}{2^n}\biggr)+\sqrt {C}
\widebar{Z}_i+\widebar{Y}_{s'}^{i,n}
\bigl(x'\bigr)\biggr)\biggr]
\end{eqnarray*}
and for $i \neq j$, $x \in I_{i,n}$, $x' \in I_{j,n}$ and $s,s' \in[n,n+1]$:
\begin{eqnarray*}
&& \mathbb E \biggl[ \biggl(X_{i,n}+\sqrt{1-\frac{1}{2^n}+C}Z+Y_{s}^{i,n}(x)
\biggr) \biggl(X_{j,n}+\sqrt {1-\frac{1}{2^n}+C}Z+Y_{s'}^{ i,n}
\bigl(x'\bigr)\biggr)\biggr]
\\[-1pt]
&&\qquad \geq \mathbb E [X_{i,n} X_{j,n}] +1- \frac{1}{2^n}+C
\\[-1pt]
&&\qquad \geq \mathbb E \biggl[\widebar{X}_{n}\biggl(\frac{i}{2^n}
\biggr) \widebar{X}_{n}\biggl(\frac
{j}{2^n}\biggr)\biggr]
\\[-1pt]
&&\qquad =\mathbb E \biggl[\biggl(\widebar{X}_{n}\biggl(
\frac{i}{2^n}\biggr)+\sqrt{C} \widebar {Z}_i+
\widebar{Y}_{s}^{i,n}(x)\biggr) \biggl(
\widebar{X}_{n}\biggl(\frac
{j}{2^n}\biggr)+\sqrt {C}
\widebar{Z}_j+\widebar{Y}_{s'}^{j,n}
\bigl(x'\bigr)\biggr)\biggr].
\end{eqnarray*}

Now, let $\beta>1$ and $r<1$ be such that $\beta r <1$ and $(\frac
{3}{2}-a)\beta r >1$. We have
{\fontsize{10pt}{12}\selectfont\begin{eqnarray*}
\hspace*{-5pt}&& \mathbb P \biggl( \sup_{0 \leq i \leq2^n-1} \sup_{s \in[n,n+1]}
\sup_{x \in
I_{i,n}} \biggl(\sqrt{2} \widebar{X}_{n}\biggl(
\frac{i}{2^n}\biggr)+\sqrt{2C} \widebar{Z}_i+\sqrt{2}\,\widebar{Y}_{s}^{i,n}(x)
\\
\hspace*{-5pt}&&\hspace*{177pt}{} -2 n \ln2 +a \ln(n+1) \biggr)\geq1\biggr)
\\
\hspace*{-5pt}&&\quad =\mathbb P \biggl( \sup_{0 \leq i \leq2^n-1} \biggl(\sqrt{2} \widebar
{X}_{n}\biggl(\frac
{i}{2^n}\biggr)+\sqrt{2C}
\widebar{Z}_i+\sqrt{2}\sup_{s \in[n,n+1],x
\in
I_{i,n}}
\widebar{Y}_{s}^{i,n}(x)
\\
\hspace*{-5pt}&&\hspace*{210pt}{} -2 n \ln2+a \ln(n+1)\biggr) \geq1
\biggr)
\\
\hspace*{-5pt}&&\quad \leq(n+1)^{a \beta r}e^{-\beta r}
\\
\hspace*{-5pt}&&\qquad{}\times   \mathbb E \Biggl[ \Biggl( \sum
_{i=0}^{2^n-1} e^{\beta
( \sqrt{2}   \widebar{X}_{n}(i/2^n)+\sqrt{2C}
\widebar{Z}_i+ \sqrt{2} \sup_{s \in[n,n+1], x \in I_{i,n}}\widebar
{Y}_{s}^{i,n}(x) -2 n \ln2 )}
\Biggr)^r \Biggr]
\\
\hspace*{-5pt}&&\quad \leq(n+1)^{a \beta r} e^{-\beta r}
\\
\hspace*{-5pt}&&\qquad{}\times  \mathbb E \Biggl[ \mathbb E
\Biggl[ \Biggl( \sum_{i=0}^{2^n-1}
e^{\beta( \sqrt{2}   \widebar{X}_{n}(i/2^n)+\sqrt{2C}
\widebar{Z}_i+ \sqrt{2} \sup_{s \in[n,n+1], x \in
I_{i,n}}\widebar
{Y}_{s}^{i,n}(x) -2 n \ln2 )} \Biggr)^r \Bigg| \widebar{X}_{n} \Biggr]
\Biggr]
\\
\hspace*{-5pt}&&\quad \leq(n+1)^{a \beta r}e^{-\beta r}
\\
\hspace*{-5pt}&&\qquad{}\times  \mathbb E \Biggl[ \Biggl( \mathbb
E \Biggl[ \sum_{i=0}^{2^n-1}
e^{\beta( \sqrt{2}   \widebar{X}_{n}(i/2^n)+\sqrt{2C}
\widebar{Z}_i+ \sqrt{2} \sup_{s \in[n,n+1], x \in
I_{i,n}}\widebar
{Y}_{s}^{i,n}(x) -2 n \ln2 )} \Bigg| \widebar{X}_{n} \Biggr] \Biggr)^r
\Biggr]
\\
\hspace*{-5pt}&&\quad \leq(n+1)^{a \beta r} e^{-\beta r} \mathbb E \bigl[ e^{\beta(\sqrt{2C}
\widebar{Z}_i+ \sqrt{2} \sup_{s \in[n,n+1], x \in
I_{i,n}}\widebar
{Y}_{s}^{i,n}(x))}
\bigr]^r
\\
\hspace*{-5pt}&&\qquad{} \times \mathbb E \Biggl[ \Biggl( \sum_{i=0}^{2^n-1}
e^{\beta( \sqrt{2}
\widebar{X}_{n}(i/2^n) -2 n \ln2 )} \Biggr)^r \Biggr]
\\
\hspace*{-5pt}&&\quad \leq C_{\beta,r} (n+1)^{a \beta r} \mathbb E \Biggl[ \Biggl( \sum
_{i=0}^{2^n-1} e^{\beta(\sqrt{2}   \widebar{X}_{n}(i/2^n) -2 n \ln2 )}
\Biggr)^r \Biggr]
\\
\hspace*{-5pt}&&\quad \leq\frac{C_{\beta,r}}{n^{(3/2-a) \beta r+o(1)}},
\end{eqnarray*}}%
%
%
where in the last line we have used Theorem 1.6 in \cite{HuShi}. This
entails the desired result by the Borel--Cantelli lemma.
\end{pf}

\subsection{Reminder about the cone construction}\label{cone}
The cone construction is based on Gaussian independently scattered
random measures; see \cite{rosinski} for further details. We consider a
Gaussian independently scattered random measure $\mu$ distributed on
the measurable space $(\mathbb{R}\times\mathbb{R}_+,\mathcal
{B}(\mathbb{R}\times\mathbb{R}_+))$, that
is, a collection of Gaussian random variables $(\mu(A),A\in\mathcal
{B}(\mathbb{R}\times\mathbb{R}_+))$ such that:
\begin{longlist}[(2)]
\item[(1)] For every sequence of disjoint sets $(A_n)_n$ in $\mathcal
{B}(\mathbb{R}
\times\mathbb{R}_+)$, the random variables $(\mu(A_n))_n$ are
independent and
\[
\mu \biggl(\bigcup_nA_n \biggr)=
\sum_n\mu(A_n)\qquad\mbox{a.s.}
\]

\item[(2)] For any measurable set $A$ in $\mathcal{B}(\mathbb{R}\times
\mathbb{R}_+)$, $\mu
(A)$ is a Gaussian random variable whose characteristic function is
given by
\[
\mathbb E \bigl(e^{iq\mu(A)}\bigr)=e^{-(q^2/2)\Gamma(A)},
\]
where the control measure $\Gamma$ is given by
\[
\Gamma(dx,dy)=\frac{1}{y^2}\,dx\,dy.
\]
We can then define the stationary Gaussian process $(\omega
_l(x))_{x\in
\mathbb{R}}$ for $0<l\leq1$ by
\[
\omega_l(x) = \mu \bigl(\mathcal{A}_l(x) \bigr),
\]
where $\mathcal{A}_l(x)$ is the triangle like subset $\mathcal
{A}_l(x):=\{(u,y)\in\mathbb{R}\times\mathbb{R}^*_+\dvtx  l \leq y \leq
1$, $-y/2 \leq x-u
\leq y/2\}$
%
\begin{figure}

\includegraphics{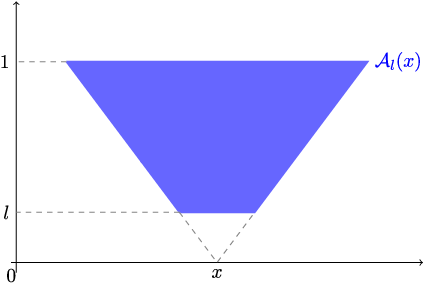}

\caption{A graphical representation of the cone construction
$\mathcal{A}_l(x)$.}\label{fig2}
%
\end{figure}
(see Figure~\ref{fig2}). The covariance kernel of the stationary Gaussian process $\omega_{l} $
is given by
%
\begin{equation}
K_l(x) = \cases{ 0, &\quad if $|x| \geq1$,
\vspace*{2pt}\cr
\displaystyle \ln
\frac{1}{|x|} + |x| - 1, &\quad if $l \leq|x|\leq 1$,
\vspace*{5pt}\cr
\displaystyle \ln
\frac{1}{l} + |x| - \frac{|x|}{l}, &\quad if $|r|\leq l$,}
\end{equation}
which can also be rewritten as
\[
K_l(x)=\int_1^{1/l}
\frac{(1-|xu|)_+}{u} \,du.
\]
Therefore, the process $\omega_{e^{-t}}$ has the same law as $X_t$. This
approach is called the cone construction.

Now we explain how to use the cone construction to prove Lemma \ref
{lemcone}, that, is to decompose the process $X_{s\ln2}=\omega
_{2^{-s}}$ for $s\in[n,n+1]$. So we choose $i\in\mathbb{N}$ such
that $0\leq
i\leq2^n-1$. We call $\mathcal{A}_{i,n}$ the common part to all the
cone like subsets $\mathcal{A}_{2^{-s}}(x)$ for $s\in[n,n+1]$ (see Figure~\ref{fig3}) and
$x\in
I_{i,n}$,
\begin{eqnarray*}
\mathcal{A}_{i,n}&=&\bigcap_{s\in[n,n+1]}\bigcap
_{x\in
I_{i,n}}\mathcal {A}_{2^{-s}}(x)
\\
&=&\biggl\{(u,y)\in\mathbb{R}\times\mathbb{R}^*_+\dvtx  2^{-n} \leq y \leq1,
-\frac{y}{2}+\frac{i+1}{2^n}\leq u\leq\frac{y}{2}+
\frac{i}{2^n }\biggr\}.
\end{eqnarray*}
For $s\in[n,n+1]$ and $x\in I_{i,n}$, we define the set $\mathcal
{R}^{i,n}_s(x)$ as
\[
\mathcal{R}^{i,n}_s(x)= \mathcal{A}_{2^{-s}}(x)
\setminus\mathcal{A}_{i,n}.
\]
\end{longlist}

\begin{figure}

\includegraphics{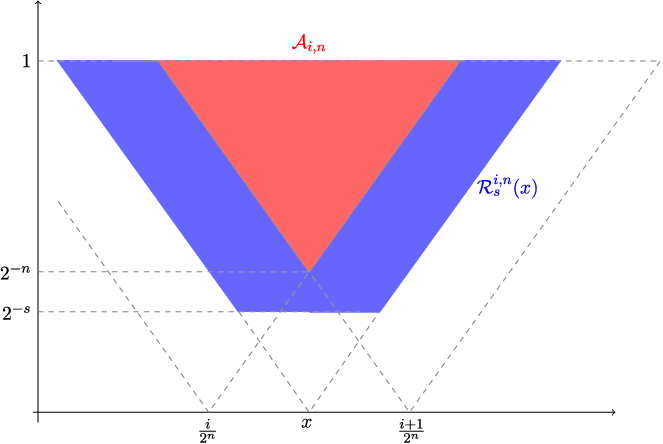}

\caption{A graphical representation of $\mathcal{A}_{i,n}$.}\label{fig3}
\end{figure}
Then we set $Y_{s}^{i,n}(x)=\mu(\mathcal{R}^{i,n}_s(x))$ and
$X_{i,n}=\mu(\mathcal{A}_{i,n})$. In particular, we find
\[
E[X_{i,n}X_{j,n}]=n\ln2+\ln\frac{1}{|i-j|+1}+
\frac{|i-j|+1}{2^{n}}-1.
\]
It is then straightforward to check the claims of Lemma~\ref{lemcone}
by using the properties of the measure $\mu$. The process
$(Y^{i,n}(x))_{s\in[n,n+1],x\in I_{i,n}}$ is independent of $X_{i,n}$
since the sets $(\mathcal{R}^{i,n}_s(x))_{s\in[n,n+1],x\in I_{i,n}}$
are all disjoint of the triangle $\mathcal{A}_{i,n}$. We also have
\[
E\bigl[Y^{i,n}(x)Y^{j,n}\bigl(x'\bigr)\bigr]
\geq0
\]
since this covariance is just given by the $\Gamma$-measure of the set
$\mathcal{R}^{i,n}_s(x)\cap\mathcal{R}^{i,n}_s(x')$. The same argument
holds to prove $E[Y^{i,n}(x)X_{j,n}]\geq0$.

\end{appendix}

\section*{Acknowledgments}
The authors wish to thank J.~Barral, J.~P.~Bouchaud, K.~Gawedzki, I.~R.~Klebanov,
I.~K.~Kostov and O.~Zindy for fruitful discussions and
comments that have led to the final version of this manuscript. The
authors also wish to thank the anonymous referee for useful comments to
improve the presentation of the paper.




\printaddresses

\end{document}